\documentclass[11pt,a4paper]{article}
\usepackage[german,english]{babel}
\usepackage{latexsym}
\usepackage{amsmath, amsfonts}
\usepackage{amssymb}	
\usepackage{amsthm}	
\usepackage{ifthen}
\usepackage{hyperref} 
\usepackage{graphicx} 

\author{Michael Klotz\footnote{Technische Universit\"at Darmstadt, Schlo\ss gartenstra\ss e 7, D-64289 Darmstadt, Deutschland,\newline klotz@mathematik.tu-darmstadt.de}}
\date{}
\title{The Automorphism Group of a Banach Principal Bundle with $\{\1\}$-structure}

\newenvironment{definition}{\begin{Definition}}{\end{Definition}\par}	
	\theoremstyle{definition} 
	\newtheorem{Definition}{Definition}[section]
\newenvironment{theorem}{\begin{Theorem}}{\end{Theorem}\par}
	\theoremstyle{plain} 
	\newtheorem{Theorem}[Definition]{Theorem}	
\newenvironment{proposition}{\begin{Proposition}}{\end{Proposition}\par}
	\newtheorem{Proposition}[Definition]{Proposition}	
\newenvironment{lemma}{\begin{Lemma}}{\end{Lemma}\par}
	\newtheorem{Lemma}[Definition]{Lemma}
\newenvironment{corollary}{\begin{Corollary}}{\end{Corollary}\par}
	\newtheorem{Corollary}[Definition]{Corollary}	

	\theoremstyle{definition}
	\newtheorem{Example}[Definition]{Example}
\newenvironment{remark}{\begin{Remark}}{\end{Remark}\par}
	\theoremstyle{definition}
	\newtheorem{Remark}[Definition]{Remark}

	\theoremstyle{definition}
	\newtheorem{Problem}[Definition]{Problem}
\renewenvironment{proof}{\noindent\textbf {Proof:}}{\hspace*{\fill} $\Box$\par\vspace{0.7ex}}

\newenvironment{acknowledgements}{\section*{Acknowledgements}}{}

\newcommand{\NN}{\mathbb N}
\newcommand{\RR}{\mathbb R}

\newcommand{\1}{{\bf 1}}

\newcommand{\calD}{{\cal D}}
\newcommand{\calL}{{\cal L}}
\newcommand{\lied}{{\cal L}} 
\newcommand{\calV}{{\cal V}}

\DeclareMathOperator{\Aut}{Aut}
\DeclareMathOperator{\Diff}{Diff}
\DeclareMathOperator{\ev}{ev}
\DeclareMathOperator{\Fix}{Fix}
\DeclareMathOperator{\flow}{Fl}
\DeclareMathOperator{\Fr}{Fr}
\DeclareMathOperator{\GL}{GL}
\DeclareMathOperator{\id}{id}
\DeclareMathOperator{\im}{im}
\DeclareMathOperator{\Iso}{Iso}
\DeclareMathOperator{\Kill}{Kill}

\newcommand{\gl}{\mathfrak{gl}}

\renewcommand{\tau}{\sigma}

\begin{document}
\maketitle

\begin{abstract}
	A $\{\1\}$-structure on a Banach manifold $M$ (with model space $E$) is an $E$-valued 1-form on $M$ that induces on each tangent space an isomorphism onto $E$. Given a Banach principal bundle $P$ with connected base space and a $\{\1\}$-structure on $P$, we show that its automorphism group can be turned into a Banach--Lie group acting smoothly on $P$ provided the Lie algebra of infinitesimal automorphisms consists of complete vector fields.
	As a consequence we show that the automorphism group of a connected geodesically complete affine Banach manifold $M$ can be turned into a Banach--Lie group acting smoothly on $M$.
	
	\noindent Keywords: Banach principal bundle, $G$-structure, automorphism group, affine Banach manifold
	
	\noindent MSC2010: 53C10, 53B05, 22F50, 22E65
\end{abstract}
%

%
%
\section{Introduction}
Given a geometric structure on a smooth manifold $M$, one of the basic associated questions is if the automorphism group can be turned into a Lie group acting smoothly on $M$. The concept of $G$-structures enables us to treat many interesting geometric structures in a unified manner. A $G$-structure on a manifold is a smooth subbundle of the bundle of linear frames with structure group $G$, where $G$ is a Lie subgroup of the general linear group of the modelling space.

The $\{\1\}$-structures on a manifold $M$ (with model space $E$) are in natural one-to-one correspondence with the fields of linear frames over $M$. Therefore, they can be defined as $E$-valued 1-forms on $M$ that induce on each tangent space an isomorphism onto $E$.

In the case of finite-dimensional manifolds, a basic theorem for $\{\1\}$-structures states that the automorphism group of a connected manifold $M$ with $\{\1\}$-structure can be turned into a Lie group acting smoothly on $M$ such that every orbit map is an injective immersion with closed image (cf.\ \cite[Th.~I.3.2]{Kob72}). 

Applying this theorem to manifolds with affine connections leads to the assertion that the automorphism group of a connected affine manifold $(M,\nabla)$ can be turned into a Lie group acting smoothly on $M$. In this situation, for each linear frame $p\in\Fr(M)$, the map $\Aut(M,\nabla)\rightarrow \Fr(M),\ f\mapsto \Fr(f)(p)$ is an injective immersion with closed image, where $\Fr(f)$ denotes the induced automorphism of the linear frame bundle $\Fr(M)$ (cf.\ \cite[Th.~II.1.3]{Kob72}).
The background to this result is that the group of affine automorphisms of $M$ is naturally isomorphic to the group of automorphisms of $\Fr(M)$ leaving both the soldering from $\theta$ and the connection form $\omega$ invariant. Note that these differential forms provide a $\{\1\}$-structure $\kappa:=(\theta,\omega)$ on the frame bundle.

A closer inspection of the proof shows that it is necessary to study $\{\1\}$-structures actually not only on connected manifolds, but, more generally, on principal bundles with connected base space, since the linear frame bundle over a connected manifold is in general not connected.

The purpose of this paper is to extend these results to smooth Banach manifolds.
This subject has also been studied by Bogdan Popescu.\footnote{
He gave a talk about Banach--Lie transformation groups at the Technical University of Darmstadt in April 2008.}
As the finite-dimensional results depend on R.~Palais' Integrability Theorem (cf.\ \cite{Pal57}), we need a comparable theorem for the infinite-dimensional case. This is provided by A.~Abouqateb and K.-H.~Neeb in \cite{AN08}. It treats the integration of locally exponential Lie algebras of complete vector fields.

Given a Banach principal bundle $P$ with connected base space and a $\{\1\}$-structure $\kappa$ on $P$, we show that the set $\Kill(P,\kappa)$ of infinitesimal automorphisms is a Banach--Lie algebra that can be embedded as a closed subspace of the tangent space at each point of $P$ by the corresponding evaluation map. If it consists of complete vector fields, it should be the natural Lie algebra of the automorphism group $\Aut(P,\kappa)$. In fact, we show that under this assumption, $\Aut(P,\kappa)$ can be turned into a Banach--Lie group acting smoothly on $P$ such that each orbit map is an injective local topological embedding.

When studying the automorphism group of a connected affine Banach manifold $(M,\nabla)$, we see that the set $\Kill(M,\nabla)$ of infinitesimal affine automorphisms inherits the structure of a Banach--Lie algebra by the given one of $\Kill(\Fr(M),\kappa)$ with $\kappa:=(\theta,\omega)$. We show that the vector fields in $\Kill(M,\nabla)$ and $\Kill(\Fr(M),\kappa)$ are complete if $M$ is assumed to be geodesically complete. The proof of these statements essentially leans on the finite-dimensional case in \cite{KN63}.
Applying our theorem about the automorphism group of a principal bundle with $\{\1\}$-structure, we show that the automorphism group of a connected affine Banach manifold $M$ that is geodesically complete can be turned into a Banach--Lie group acting smoothly\linebreak on $M$.

The original motivation of this research is to attack a problem concerning symmetric spaces in the sense of O.~Loos (cf.\ \cite{Loo69}) that are modelled on Banach spaces.
We can show that the automorphism group $G$ of a connected symmetric space $M$ is a Banach--Lie group acting transitively on $M$. In particular, $M$ is a Banach homogeneous space. More precisely, we have $M\cong G/G_b$, where the stabilizer $G_b$ for a point $b\in M$ is an open subgroup of the group of fixed points in $G$ for the involution $\sigma$ on $G$ given by $\sigma(g):=\mu_b\circ g\circ \mu_b$ with the symmetry $\mu_b$ at $b$ (cf.\ \cite[Ex.~3.9]{Nee02} for homogeneous symmetric spaces).
As a connected symmetric space carries a canonical affine connection encoding the symmetric structure in the sense that it has the same automorphisms, we can apply the results of this paper to show its homogeneity. For details we refer to \cite{Klo09_Sym}.
\tableofcontents
%
%
%
%
\section{The Automorphism Group of a Banach Principal Bundle with \{1\}-structure}
\subsection{Introduction}\label{sec:subsectionIntro}
Let $(P,M,G,q,\rho)$ be a smooth Banach principal bundle, i.e., $q\colon P\rightarrow M$ is a smooth map of Banach manifolds and $\rho\colon P\times G \rightarrow P$ a smooth action of a Banach--Lie group $G$ on $P$ with the property of local triviality: Each $x\in M$ has an open neighborhood $U$ for which there exists a smooth diffeomorphism $\varphi\colon q^{-1}(U)\rightarrow U \times G$ satisfying
$$q(\varphi^{-1}(x,g))=x \quad\mbox{and}\quad \varphi^{-1}(x,g_1g_2)=\varphi^{-1}(x,g_1).g_2$$
for all $x\in U$ and $g, g_1,g_2\in G$ (cf.\ \cite[6.2.1]{Bou07}).

We assume $M$ to be pure with model space $E$ and denote the Lie algebra of $G$ by $\mathfrak{g}$.
A \emph{$\{\1\}$-structure} on $P$ is an $(E\times\mathfrak{g})$-valued 1-form
$\kappa$
on $P$ such that $\kappa_p\colon T_pP\rightarrow E\times \mathfrak{g}$ is a topological linear isomorphism for each $p\in P$.

We denote by $\Aut(P,\kappa)$ the group of principal bundle automorphisms of $P$ that leave $\kappa$ invariant, i.e.,
$$\Aut(P,\kappa)\ :=\ \{f\in\Diff(P)\colon f^\ast\kappa=\kappa \mbox{ and } f\circ\rho_g=\rho_g\circ f \mbox{ for all } g\in G\}.$$
The set $\Kill(P,\kappa)$ of \emph{infinitesimal automorphisms} is defined by
$$\Kill(P,\kappa) \ :=\ \{\xi\in\calV(P)\colon \lied_\xi\kappa=0 \mbox{ and } (\rho_g)_\ast\xi=\xi \mbox{ for all } g\in G\}.$$
In the light of Section~\ref{sec:basicConcepts}, this is the set of all smooth vector fields $\xi\in\calV(P)$ whose flow maps $\flow^\xi_t\colon\calD_t(\xi)\rightarrow\calD_{-t}(\xi)$ leave $\kappa$ invariant and commute with $\rho_g$ for all $g\in G$.

Our main results are:
\begin{proposition}\label{prop:Banach-LieAlgebraKill(M,kappa,rho)}
	Let $(P,M,G,q,\rho)$ be a smooth Banach principal bundle and $\kappa$ a $\{\1\}$-structure on $P$. The set $\Kill(P,\kappa)$ of infinitesimal automorphisms is a Lie subalgebra of the Lie algebra $\calV(P)$ of smooth vector fields. If $M$ is connected, then it carries a unique Banach space structure such that each evaluation map $\ev_p\colon \Kill(P,\kappa)\rightarrow T_pP,\ \xi\mapsto\xi(p)$ is a closed embedding. With this structure it becomes a Banach--Lie algebra.
\end{proposition}
\begin{theorem}\label{th:lieGroupKill(M,kappa,rho)}
	Let $(P,M,G,q,\rho)$ be a smooth Banach principal bundle, where the base space $M$ is connected, and let $\kappa$ be a $\{\1\}$-structure on $P$. If all infinitesimal automorphisms $\xi\in\Kill(P,\kappa)$ are complete, then the automorphism group $\Aut(P,\kappa)$ can be turned into a Banach--Lie group such that
	$$\exp\colon \Kill(P,\kappa)\rightarrow \Aut(P,\kappa),\ \xi \mapsto \flow^{-\xi}_1$$
 	is its exponential map. The natural map $\tau\colon \Aut(P,\kappa)\times P \rightarrow P$ is a smooth action whose derived action is the inclusion map $\Kill(P,\kappa)\hookrightarrow\calV(P)$, i.e., $-T\tau(\id_P,p)(\xi,0)=\xi(p)$.
 	For each $p\in P$, the orbit map $\tau_p\colon \Aut(P,\kappa) \rightarrow P,\ g\mapsto g(p)$ is an injective local topological embedding.
\end{theorem}
When regarding the trivial bundle $P=M\times G$ with trivial structure group $G=\{\1\}$, we obtain the following corollary:
\begin{corollary}
	Let $M$ be a connected smooth Banach manifold and $\kappa$ a $\{\1\}$-structure on $M$. If all infinitesimal automorphisms $\xi\in\Kill(M,\kappa)$ are complete, then the automorphism group $\Aut(M,\kappa)$ can be turned into a Banach--Lie group acting smoothly on $M$.
\end{corollary}
One consequence (of the theorem) that we shall show is that the automorphism group of a geodesically complete affine Banach manifold $M$ can be turned into a Banach--Lie group that acts smoothly on $M$. 
\begin{remark}
	In \cite{Kob72}, S.~Kobayashi deals with the finite-dimensional case. He shows that the automorphism group of a connected manifold $M$ with a $\{1\}$-structure can be turned into a Lie group acting smoothly on $M$ such that each orbit map is an injective immersion with closed image (cf.\ \cite[Th.~I.3.2]{Kob72}).
\end{remark}
\subsection{Basic Concepts} \label{sec:basicConcepts}
\begin{lemma}[{cf.\ \cite[Lem~3.3]{AN09}}] \label{lem:T_xf=closedEmbedding}
	Let $M$ and $N$ be smooth Banach manifolds, $U$ an open neighborhood of some point $x$ in $M$ and $f\colon U\rightarrow N$ a $C^1$-map for which $T_xf$ is a closed embedding. Then there exists an open $x$-neighborhood $U^\prime\subseteq U$ such that $f|_{U^\prime}$ is a topological embedding.
\end{lemma}
\begin{definition}[Lie derivative]
	Let $M$ be a smooth Banach manifold and $\xi\in\calV(M)$ a smooth vector field on $M$ with local flow $\flow^\xi\colon\RR\times M\supseteq \calD(\xi)\rightarrow M$. Further, let $\omega$ be a smooth vector field on $M$ or a Banach space valued differential form on $M$ of any degree. The \emph{Lie derivative $\lied_\xi\omega$ of $\omega$ with respect to $\xi$} is the smooth vector field or the differential form, respectively, given by
	$$(\lied_\xi\omega)(x) \ =\ \left.\frac{d}{dt}\right|_{t=0}\big((\flow^\xi_t)^\ast\omega\big)(x) \ =\ \lim_{t\rightarrow 0}\frac{\big((\flow^\xi_t)^\ast\omega\big)(x)-\omega(x)}{t}$$
\end{definition}
For $\eta\in\calV(M)$, we have $\lied_\xi\eta=[\xi,\eta]$, i.e., the Lie bracket in $\calV(M)$, and for a differential form $f$ of degree 0, i.e., for a smooth function $f$ on $M$, we have $\lied_\xi f=\xi.f:=df\circ \xi$ (cf.\ \cite[V, \S 2]{Lan01}). Given a differential form $\omega$ of degree $n\geq 1$ and $\eta_1,\dots,\eta_n\in\calV(M)$, we have
$$\lied_\xi(\omega(\eta_1,\dots, \eta_n)) \ =\ (\lied_\xi\omega)(\eta_1,\dots,\eta_n)+\sum_{i=1}^n\omega(\eta_1,\dots,\lied_\xi\eta_i,\dots,\eta_n)$$
(cf.\ \cite[Prop.~V.5.1]{Lan01}).

\begin{proposition}[{cf.\ \cite[Prop.~4.2.4]{AMR88}}]\label{prop:fRelatedVectorFields}
	Let $M$ and $N$ be smooth Banach manifolds, $\xi\in\calV(M)$ and $\eta\in\calV(N)$ smooth vector fields and $f\colon U\rightarrow V$ a smooth map between open submanifolds $U\subseteq M$ and $V\subseteq N$. Then the following conditions are equivalent:
	\begin{enumerate}
		\item[\rm (a)] $\xi|_U$ and $\eta|_V$ are $f$-related, i.e., $Tf\circ\xi|_U=\eta\circ f$, also denoted by $f_\ast(\xi|_U)=\eta|_V$.
		\item[\rm (b)] $f$ intertwines the flow maps $\flow^\xi_t$ and $\flow^\eta_t$  in the sense that
		$$(f\circ\flow^\xi_t)(x) \ =\ (\flow^\eta_t\circ f)(x)$$
		 for all $(t,x)\in\calD(\xi)$ that satisfy $\flow^\xi([0,t],x)\subseteq U$ or $\flow^\xi([t,0],x)\subseteq U$, respectively.
	\end{enumerate}
	If we have $U=M$ and $V=N$, then {\rm (b)} means
		$$f(\calD_t(\xi))\subseteq\calD_t(\eta) \quad\mbox{and}\quad f\circ\flow^\xi_t=\flow^\eta_t\circ f|_{\calD_t(\xi)}.$$
\end{proposition}
%
%
\begin{proposition}[{cf.\ \cite[Cor.~5.4.2]{AMR88}}]\label{prop:lieXiOmega=0}
	Let $M$ be a smooth Banach manifold and $\xi\in\calV(M)$ a smooth vector field. Further, let $\omega$ be a smooth vector field on $M$ or a Banach space valued differential form on $M$ of any degree. The following conditions are equivalent:
	\begin{enumerate}
		\item[\rm (a)] $\lied_\xi\omega=0.$
		\item[\rm (b)] $\omega$ is invariant under all flow maps $\flow^\xi_t\colon\calD_t(\eta)\rightarrow\calD_{-t}(\eta)$, i.e., $(\flow^\xi_t)^\ast(\omega|_{\calD_{-t}(\xi)})=\omega|_{\calD_t(\xi)}$.
	\end{enumerate}
\end{proposition}
%
%
\begin{corollary} \label{cor:commutingFlows}
	Let $M$ be a smooth Banach manifold and $\xi,\eta\in\calV(M)$ smooth vector fields. The following conditions are equivalent:
	\begin{enumerate}
		\item[\rm (a)] $[\xi,\eta]=0.$
		\item[\rm (b)] $\eta$ is invariant under all flow maps $\flow^\xi_t\colon\calD_t(\eta)\rightarrow\calD_{-t}(\eta)$, i.e., $(\flow^\xi_t)_\ast\eta|_{\calD_t(\xi)}=\eta|_{\calD_{-t}(\xi)}$.
		\item[\rm (c)] The flow maps $\flow^\eta_t$ and $\flow^\xi_s$ commute in the sense that
		$$(\flow^\xi_s\circ\flow^\eta_t)(x) \ =\ (\flow^\eta_t\circ \flow^\xi_s)(x)$$
		 for all $(t,x)\in\calD(\eta)$ and $s\in\RR$ that satisfy $\flow^\eta([0,t],x)\subseteq \calD_s(\xi)$ or $\flow^\eta([t,0],x)\subseteq \calD_s(\xi)$, respectively.
	\end{enumerate}
\end{corollary}
\begin{lemma} \label{lem:parameterDependingFlow}
	Let $\xi\colon F \times M \rightarrow TM$ be a smooth vector field depending on parameters in a Banach space $F$ such that $F\rightarrow\calV(M),\ v\mapsto \xi_v:=\xi(v,\cdot)$ is linear. Let $\flow^\xi\colon\RR\times F\times M \supseteq\calD(\xi)\rightarrow M$ be its local flow. Then, for each $x\in M$, the derivative of the map $\flow^\xi_{1,x}:=\flow^\xi(1,\cdot,x)$ at $0$ is given by $T_0\flow^\xi_{1,x}=\xi(\cdot,x)\colon F\rightarrow T_xM$.
\end{lemma}
\begin{proof}
	We have
	\begin{eqnarray*}
		(T_0\flow^\xi_{1,x})(v) &=& \left.\frac{d}{dt}\right|_{t=0}\flow^\xi_{1,x}(0+tv) \ =\ \left.\frac{d}{dt}\right|_{t=0}\flow^{\xi_{tv}}_1(x)\\
		 &=& \left.\frac{d}{dt}\right|_{t=0}\flow^{t\xi_{v}}_1(x) \ =\ \left.\frac{d}{dt}\right|_{t=0}\flow^{\xi_{v}}_t(x) \ =\ \xi_v(x).
	\end{eqnarray*}
\end{proof}
\begin{theorem}[cf.\ \cite{AN08}] \label{th:integrationOfLieAlgebrasOfVectorFields}
	Let $\mathfrak{g}$ be a Banach--Lie algebra, $M$ a smooth Banach manifold and $\alpha \colon \mathfrak{g}\rightarrow \calV(M)$ an injective morphism of Lie algebras satisfying
	\begin{enumerate}
 		\item[\rm (1)] Each vector field $\alpha(X)$ is complete.
 		\item[\rm (2)] The map $\widehat\alpha:\mathfrak{g}\times M\rightarrow TM,\ (X,x)\mapsto \alpha(X)(x)$ is smooth.
 		\item[\rm (3)] The subgroup $\Gamma_\alpha:=\{X\in\mathfrak{z}(\mathfrak{g})\colon \flow^{\alpha(X)}_1=\id_M\}$ of the center $\mathfrak{z}(\mathfrak{g})$ of $\mathfrak{g}$ is discrete.
 	\end{enumerate}
 	Then the subgroup $G:=\langle\flow^{\alpha(X)}_1\colon X\in\mathfrak{g}\rangle_{\Diff(M)}$ of the diffeomorphism  group $\Diff(M)$ can be equipped with a unique Banach--Lie group structure such that
 	$$\exp_\alpha\colon \mathfrak{g}\rightarrow G,\ X \mapsto \flow^{-\alpha(X)}_1$$
 	is its exponential map. The natural map $\beta\colon G\times M \rightarrow M$ is a smooth action with $\dot{\beta}=\alpha$.
\end{theorem}
Here $\dot{\beta}\colon \mathfrak{g}\rightarrow\calV(M)$ is the derived action defined by
$$\dot{\beta}(X)(x):=-T\beta(\id_M,x)(X,0).$$
\begin{remark}
	In \cite{AN08}, condition (2) is that the map $\mathfrak{g} \times M\rightarrow M, \ (X,x)\mapsto \flow^{\alpha(X)}_1(x)$ is smooth. In the Banach case, this is a conclusion of our condition, as this map is the time-1-flow of the smooth vector field $\widehat\alpha$ that depends on parameters in $\mathfrak{g}$ (cf.\ \cite[pp.\ 72, 92, 160]{Lan01}).
\end{remark}
\subsection{The Banach--Lie Algebra of the Infinitesimal Automorphisms} \label{sec:infinitesimalAutomorphisms}
Let $(P,M,G,q,\rho)$ be a smooth Banach principal bundle where the base space $M$ (modelled on $E$) is connected and let $\kappa$ be a $\{\1\}$-structure on $P$.
For each $v\in F:= E\times \mathfrak{g}$, let $\eta_v\in\calV(P)$ be the smooth vector field defined by $$\eta_v(p):=\kappa_p^{-1}(v).$$
It is clear that $\kappa(\eta_v)(p)\equiv v$. We shall often use that, for each $p\in P$, the tangent space $T_pP$ is given by $\{\eta_v(p)\colon v\in F\}$.
The following lemma is also valid for non-connected $M$.
\begin{lemma}\label{lem:Kill(M,kappa,rho)}
	We have
	$$\Kill(P,\kappa) \ =\ \{\xi\in\calV(P)\colon [\xi,\eta_v]=0 \mbox{ and } (\rho_g)_\ast\xi=\xi \mbox{ for all } v\in F \mbox{ and } g\in G\}.$$
	It follows that $\Kill(P,\kappa)$ is a Lie subalgebra of the Lie algebra $\calV(P)$.
\end{lemma}
\begin{proof}
	For any $\xi\in\calV(P)$, we have $\lied_\xi\kappa=0$ if and only if $(\lied_\xi\kappa)(\eta_v)=0$ for all $v\in F$. Due to $(\lied_\xi\kappa)(\eta_v)=\xi.(\kappa(\eta_v))-\kappa([\xi,\eta_v])= 0 -\kappa([\xi,\eta_v])$, this is equivalent to $[\xi,\eta_v]=0$ for all $v\in F$, the map $\kappa_p$ being injective for all $p\in P$. We thus obtain the first assertion.
	An easy computation using the Jacobi identity and the naturality of the Lie bracket shows that $\Kill(P,\kappa)$ is a Lie subalgebra of $\calV(P)$.
\end{proof}
Let $\eta\colon F\times P\rightarrow TP$ be the smooth vector field depending on parameters that is given by
$$\eta(v,p):= \eta_v(p)=\kappa_p^{-1}(v)$$
and let $\flow^\eta\colon\RR\times F\times P\supseteq \calD(\eta)\rightarrow P$ be its local flow. The map $F\rightarrow \calV(P),\ v\mapsto \eta_v$ being linear, Lemma~\ref{lem:parameterDependingFlow} says that, for each $p\in P$, the derivative of the map $\flow^\eta_{1,p}:=\flow^\eta(1,\cdot,p)$ at 0 is given by $T_0\flow^\eta_{1,p}=\eta(\cdot,p)=\kappa_p^{-1}$, i.e., by a topological linear isomorphism. Hence, the map $\flow^\eta_{1,p}$ induces a local diffeomorphism at $0\in F$.

Therefore, for each $p\in P$, there is a diffeomorphism
\begin{equation} \label{eqn:phiX(v)}
	\varphi_p\colon V_p\rightarrow U_p\subseteq P\quad\mbox{with}\quad \varphi_p(v)=\flow^{\eta_v}_1(p)\quad\mbox{and}\quad T_0\varphi_p=\kappa_p^{-1}
\end{equation}
of an open neighborhood $V_p$ of $0\in F$ onto an open neighborhood $U_p$ of $p\in P$.

We define the set
\begin{equation} \label{eqn:setS}
	S:=\big\{s:=s_n\circ\cdots\circ s_1\colon n\in\NN \mbox{ and } s_1,\dots, s_n \in \{\flow^{\eta_v}_t\colon t\in\RR, v\in F\}\cup\{\rho_g\colon g\in G\} \big\},
\end{equation}
where each $s_n\circ\cdots\circ s_1$ means the composition with maximal domain.
It is stable under composition and inversion of maps.
\begin{lemma} \label{lem:s*Xi=Xi,sInS}
	Given $\xi\in\calV(P)$, we have $\xi\in \Kill(P,\kappa)$ if and only if for each  $s\colon U\rightarrow V$ in $S$, we have $s_\ast(\xi|_U)=\xi|_V$, i.e., $\xi(s(p))=Ts(\xi(p))$ for all $p\in U$.
\end{lemma}
\begin{proof}
	By Lemma~\ref{lem:Kill(M,kappa,rho)} and Corollary~\ref{cor:commutingFlows}, we have $\xi\in \Kill(P,\kappa)$ if and only if $(\rho_g)_\ast\xi=\xi$ and $(\flow^{\eta_v}_t)_\ast(\xi|_{\calD_t(\eta_v)})=\xi_{\calD_{-t}(\eta_v)}$ for all $g\in G$, $t\in\RR$ and $v\in F$.
	As the set $S$ consists of compositions of such maps, the assertion follows.
\end{proof}
\begin{lemma}\label{lem:SconnectsM}
	Given any $p_1, p_2\in P$, there is a map $s\in S$ with $s(p_1)=p_2$.
\end{lemma}
\begin{proof}
	For a fixed $p_1\in P$, we put $A:=\{p\in P\colon (\exists s\in S)(s(p_1)=p)\}$. We shall show that $A$ is all of $P$. It is clear that $A$ is $\rho$-invariant and not empty. As the base space $M$ is connected, it suffices to show that $A$ and its complement $A^c$ both are open in $P$.
	
	For any $p\in A$ (with $s(p_1)=p$), the neighborhood $U_p$ (cf.\ (\ref{eqn:phiX(v)})) is contained in $A$, since for each $\varphi_p(v)\in U_p$, we have $\varphi_p(v)=(\flow^{\eta_v}_1\circ s)(p_1)$.
	Hence, $A$ is open. To see that $A^c$ is open, let $p$ be any point in $A^c$. We shall show that $U_p\subseteq A^c$. If any $\varphi_p(v)\in U_p$ was in $A$, i.e., $s(p_1)=\varphi_p(v)$ for some $s\in S$, then $p$ was in $A$, too, due to $p = (\flow^{\eta_{-v}}_1\circ s)(p_1)$.
\end{proof}
For each $p\in P$, let $H_p$ be the subgroup of $\GL(T_pP)$ defined by
$$H_p:=\{T_ps\colon s\in S \mbox{ and } s(p)=p\}.$$
\begin{proposition}\label{prop:Kill=BanachSpace}
	For each $p\in P$, the evaluation map $\ev_p\colon \Kill(P,\kappa)\rightarrow T_pP,\ \xi\mapsto\xi(p)$ is an injective linear map. Its image $\im(\ev_p)$ is given by the fixed point set $(T_pP)^{H_p}$ under the group $H_p\leq \GL(T_pP)$. For each $p_1,p_2\in P$, the map $\ev_{p_2}\circ(\ev_{p_1})^{-1}\colon (T_{p_1}P)^{H_{p_1}}\rightarrow (T_{p_2}P)^{H_{p_2}}$ is a topological linear isomorphism of Banach spaces. It is given by restricting $T_{p_1}s$ for any $s\in S$ satisfying $s(p_1)=p_2$.
\end{proposition}
\begin{proof}
	To see the injectivity of $\ev_p$, we shall show that for each $\xi\in\Kill(P,\kappa)$ and $p^\prime\in P$, the vector $\xi(p^\prime)$ depends only on $\xi(p)$. Indeed, by Lemma~\ref{lem:SconnectsM} and Lemma~\ref{lem:s*Xi=Xi,sInS}, there is a map $s\in S$ with $s(p)=p^\prime$ and $\xi(p^\prime)=Ts(\xi(p))$.
		
	The inclusion $\im(\ev_p)\subseteq(T_pM)^{H_p}$ is clear by Lemma~\ref{lem:s*Xi=Xi,sInS}. Conversely, for each $w\in(T_pP)^{H_p}$, let $\xi_w\in\calV(P)$ be the smooth vector field defined by
	\begin{equation}\label{eqn:defXiW}
		\xi_w(s(p)):=Ts(w)\quad\mbox{for all $s\in S$ where $s(p)$ is defined.}
	\end{equation}
	Due to Lemma~\ref{lem:SconnectsM}, $\xi_w$ is defined on all of $P$. It is well-defined, as
	$$s_1(p)=s_2(p) \ \Rightarrow\ ((s_2)^{-1}\circ s_1)(p)=p \ \Rightarrow\ T((s_2)^{-1}\circ s_1)(w)=w \ \Rightarrow\ Ts_1(w)=Ts_2(w)$$
	for all $s_1,s_2 \in S$. Before checking its smoothness, we observe that
	\begin{equation}\label{eqn:xi_w(s(y))}
	\xi_w(s(p^\prime))=Ts(\xi_w(p^\prime))\quad\mbox{for all $s\in S$ and $p^\prime\in P$ where $s(p^\prime)$ is defined},
	\end{equation}
	since
	$$\xi_w(s(p^\prime))\ =\ \xi_w((s\circ r)(p))\ =\ (Ts\circ Tr)(w)\ =\ Ts\big(\xi_w(r(p))\big)\ =\ Ts(\xi_w(p^\prime))$$
	for an appropriate $r\in S$ with $r(p)=p^\prime$.
	
	To see its smoothness, we work locally on $U_{p^\prime}$ (cf.\ (\ref{eqn:phiX(v)})) for all $p^\prime\in P$. We have
	\begin{equation} \label{eqn:xi_w(phi_y(v))}
	\xi_w(\varphi_{p^\prime}(v))\ =\ \xi_w(\flow^{\eta_v}_1(p^\prime))\ =\ T\flow^{\eta_v}_1(\xi_w(p^\prime))\ =\ (T_{(v,p^\prime)}\flow^{\eta}_1)(0,\xi_w(p^\prime))
	\end{equation}
	for all $v\in V_{p^\prime}$. Hence, $\xi_w|_{U_{p^\prime}}$ is smooth. By (\ref{eqn:xi_w(s(y))}) and Lemma~\ref{lem:s*Xi=Xi,sInS}, we thus know that $\xi_w\in\Kill(P,\kappa)$.
	Because of $\xi_w(p)=w$, we then obtain $(T_pP)^{H_p}\subseteq\im(\ev_p)$.
	
	Given any $s\in S$ with $s(p_1)=p_2$, we have $$(\ev_{p_2}\circ(\ev_{p_1})^{-1})(w)\ =\ \ev_{s(p_1)}(\xi_w)\ \stackrel{\mbox{\scriptsize(\ref{eqn:defXiW})}}{=} \ Ts(w)$$
	for all $w\in (T_pP)^{H_p}$. In particular, $\ev_{p_2}\circ(\ev_{p_1})^{-1}$ is a topological linear isomorphism.
\end{proof}
\begin{corollary} \label{cor:Banach-LieAlgebraKill(M,kappa,rho)}
	The Lie subalgebra $\Kill(P,\kappa)$ of the Lie algebra $\calV(P)$ carries a unique Banach space structure such that each evaluation map $\ev_p\colon \Kill(P,\kappa)\rightarrow T_pP,\ \xi\mapsto\xi(p)$ is a closed embedding. With this structure it becomes a Banach--Lie algebra.
\end{corollary}
\begin{proof}
	It remains to check the compatibility of the Banach space structure and the Lie bracket of $\Kill(P,\kappa)$. To see that
	$$\Kill(P,\kappa)\times\Kill(P,\kappa)\rightarrow\Kill(P,\kappa),\quad (\xi_{w_1},\xi_{w_2})\mapsto [\xi_{w_1},\xi_{w_2}]$$
	is continuous, we consider one embedding $\ev_p$ and have to check the continuity of
	$$(T_pP)^{H_p}\times T_pP^{H_p} \rightarrow T_pP^{H_p}, \quad (w_1,w_2)\rightarrow [\xi_{w_1},\xi_{w_2}](p).$$
	We shall work in the chart $\varphi_p^{-1}$ (cf.\ (\ref{eqn:phiX(v)})). For this, we restrict the flow map
	$\flow^\eta_1\colon F\times P\supseteq \calD_1(\eta)\rightarrow P$ to a set $\widetilde V_p \times \widetilde U_p$ with open neighborhoods $\widetilde V_p$ of $0\in F$ and $\widetilde U_p$ of $p\in P$ such that $\flow^\eta_1(\widetilde V_p \times \widetilde U_p)\subseteq U_p$ and $\widetilde\varphi_p^{-1}\colon \widetilde U_p\rightarrow \widetilde V_p\subseteq F$ is a restriction of the chart $\varphi_p^{-1}$. This works, because $\flow^\eta_1(0,p)=p$. The flow map $\flow^\eta_1$ then has a local representation
	$(\flow^\eta_1)^{\widetilde\varphi_p}\colon\widetilde V_p \times \widetilde V_p\rightarrow V_p$.
	
	Due to (\ref{eqn:xi_w(phi_y(v))}), the local representations $\xi_{w}^{\widetilde\varphi_p}$ of $\xi_w$ are given by
	$$\xi_{w}^{\widetilde\varphi_p}(v) \ =\ d_2(\flow^\eta_1)^{\widetilde\varphi_p}(v,0)(\kappa_p(w))$$
	with derivatives
	$$d\xi_{w}^{\widetilde\varphi_p}(v)(u) \ =\ d_1d_2(\flow^\eta_1)^{\widetilde\varphi_p}(v,0)(u)(\kappa_p(w)).$$
	We shall verify the continuity of
	$$\kappa_p\big((T_pP)^{H_p}\big)\times \kappa_p\big((T_pP)^{H_p}\big) \rightarrow \kappa_p\big((T_pP)^{H_p}\big), \quad (\bar w_1,\bar w_2)\rightarrow [\xi^{\widetilde\varphi_p}_{w_1},\xi^{\widetilde\varphi_p}_{w_2}](0),$$
	with $w_i:=\kappa_p^{-1}(\bar w_i)$, for $i=1,2$.
	We have
	\begin{eqnarray*}
		[\xi^{\widetilde\varphi_p}_{w_1},\xi^{\widetilde\varphi_p}_{w_2}](0) &=& d\xi^{\widetilde\varphi_p}_{w_2}(0)(\xi^{\widetilde\varphi_p}_{w_1}(0))-d\xi^{\widetilde\varphi_p}_{w_1}(0)(\xi^{\widetilde\varphi_p}_{w_2}(0))\\
		&=& d_1d_2(\flow^\eta_1)^{\widetilde\varphi_p}(0,0)(\bar w_1)(\bar w_2)
			- d_1d_2(\flow^\eta_1)^{\widetilde\varphi_p}(0,0)\big(\bar w_2)(\bar w_1),
	\end{eqnarray*}
	which depends smoothly on $(\bar w_1,\bar w_2)$.
\end{proof}
\subsection{Proof of the Main Theorem}
We consider the situation of Section~\ref{sec:infinitesimalAutomorphisms} and use the same notation.
\begin{lemma} \label{lem:Aut(M,kappa,rho)}
	We have
	\emph{$$\Aut(P,\kappa)\ =\ \{f\in\Diff(P)\colon f_\ast\eta_v=\eta_v \mbox{ and } f\circ\rho_g=\rho_g\circ f \mbox{ for all $v\in F$ and $g\in G$} \}.$$}
\end{lemma}
\begin{proof}
	Given $f\in\Diff(P)$, we have $f^\ast\kappa = \kappa$ if and only if
	$(f^\ast\kappa)_p(\eta_v(p))=\kappa_p(\eta_v(p))$ for all $p\in P$ and $v\in F$, i.e., $\kappa_{f(p)}(Tf(\eta_v(p)))=v$. This is equivalent to $Tf(\eta_v(p))=\kappa_{f(p)}^{-1}(v)$, i.e., $Tf(\eta_v(p))=\eta_v(f(p))$. This means $f_\ast\eta_v=\eta_v$ for all $v\in F$. The assertion follows from the definition of $\Aut(P,\kappa)$.
\end{proof}
\begin{lemma} \label{lem:fs=sf,sInS}
	Given $f\in\Diff(P)$, we have $f\in\Aut(P,\kappa)$ if and only if for each $s\colon U\rightarrow V$ in $S$ (cf.\ \emph{(\ref{eqn:setS})}), we have $f(U)\subseteq U$ and $f \circ s = s\circ f|_U$.
\end{lemma}
\begin{proof}
	By Lemma~\ref{lem:Aut(M,kappa,rho)} and Proposition~\ref{prop:fRelatedVectorFields}, we have $f\in\Aut(P,\kappa)$ if and only if $f\circ\rho_g=\rho_g\circ f$, $f(\calD_t(\eta_v))\subseteq \calD_t(\eta_v)$ and $f\circ\flow^{\eta_v}_t=\flow^{\eta_v}_t\circ f|_{\calD_t(\eta_v)}$ for all $g\in G$, $t\in \RR$ and $v\in F$. As the set $S$ consists of compositions of such maps, the assertion follows. 
	Notice that, given $s_1\colon U_1\rightarrow V_1$ and $s_2\colon U_2\rightarrow V_2$ satisfying the condition, the domain $s_1^{-1}(U_2)$ of the composition $s_2\circ s_1|_{s_1^{-1}(U_2)}$ satisfies $f(s_1^{-1}(U_2))\subseteq s_1^{-1}(U_2)$, since
	$$s_1\big(f(s_1^{-1}(U_2))\big)\ =\ f\big(s_1(s_1^{-1}(U_2))\big)\ \subseteq\ f(U_2)\ \subseteq\ U_2.$$
\end{proof}
\begin{lemma}
 	If all infinitesimal automorphisms in $\Kill(P,\kappa)$ are complete, then the inclusion map $\Kill(P,\kappa)\hookrightarrow \calV(P)$ satisfies the conditions of Theorem~\ref{th:integrationOfLieAlgebrasOfVectorFields}:
 	\begin{enumerate}
 		\item[\rm (1)] Each vector field $\xi\in\Kill(P,\kappa)\subseteq \calV(P)$ is complete (by assumption).
 		\item[\rm (2)] The evaluation map $\ev\colon\Kill(P,\kappa) \times P\rightarrow TP, \ (\xi,p)\mapsto \xi(p)$ is smooth.
 		\item[\rm (3)] The subgroup $\Gamma:=\{\xi\in\mathfrak{z}(\Kill(P,\kappa))\colon \flow^\xi_1=\id_P\}$ of the center $\mathfrak{z}(\Kill(P,\kappa))$ of $\Kill(P,\kappa)$ is discrete.
 	\end{enumerate}
\end{lemma}
\begin{proof}
	(2) Considering an embedding $\ev_p\colon \Kill(P,\kappa)\hookrightarrow T_pP$, we have to check the smoothness of
	$$(T_pP)^{H_p}\times P \rightarrow TP,\ (w,p^\prime)\rightarrow \xi_w(p^\prime),$$
	where $\xi_w:=\ev_p^{-1}(w)$ (cf.\ Proposition~\ref{prop:Kill=BanachSpace}). We work locally on $(T_pP)^{H_p}\times U_{p^\prime}$ (cf.\ (\ref{eqn:phiX(v)})) for all $p^\prime\in P$. By Lemma~\ref{lem:SconnectsM}, there is an $s\in S$ with $s(p)=p^\prime$. From (\ref{eqn:xi_w(phi_y(v))}) and (\ref{eqn:defXiW}) (in Section~\ref{sec:infinitesimalAutomorphisms}), we know that
	$$\xi_w(\varphi_{p^\prime}(v))\ =\ (T_{(v,p^\prime)}\flow^{\eta}_1)\big(0,\xi_w(s(x))\big)\ =\ (T_{(v,p^\prime)}\flow^{\eta}_1)(0,Ts(w)),$$
	which depends smoothly on $(w,v)$.
	
	(3) It suffices to show that there is a 0-neighborhood $U\subseteq\Kill(P,\kappa)$ with $\flow^\xi_1\neq\id_P$ for all $\xi\in U\backslash \{0\}$. Fixing some $p\in P$ and considering the map $\psi_p\colon\Kill(P,\kappa)\rightarrow M,\ \xi \mapsto \flow^\xi_1(p)$, we shall show that $\psi_p$ is injective on an appropriate 0-neighborhood $U$, so that $\flow^\xi_1(p)\neq \id_P(p)$ and hence $\flow^\xi_1\neq \id_P$ for all $\xi\in U\backslash \{0\}$.
	Due to (2), we can apply Lemma~\ref{lem:parameterDependingFlow} and obtain
	$(T_0\psi_p)=\ev_p\colon\Kill(P,\kappa)\rightarrow T_pP$, which is a closed embedding. From Lemma~\ref{lem:T_xf=closedEmbedding}, we then know that, for an appropriate 0-neighborhood $U$, $\psi_p|_U$ is a topological embedding and hence injective.	
\end{proof}
As a consequence of Theorem~\ref{th:integrationOfLieAlgebrasOfVectorFields}, we obtain:
\begin{corollary}
	If all infinitesimal automorphisms in $\Kill(P,\kappa)$ are complete, then the subgroup $H:=\langle\flow^{\xi}_1\colon \xi\in\Kill(P,\kappa)\rangle_{\Diff(P)}$ of the diffeomorphism  group $\Diff(P)$ can be equipped with a unique Banach--Lie group structure such that
 	$$\exp\colon \Kill(P,\kappa)\rightarrow H,\ \xi \mapsto \flow^{-\xi}_1$$
 	is its exponential map. The natural map $\beta\colon H\times P \rightarrow P$ is then a smooth action whose derived action is the inclusion map $\Kill(P,\kappa)\hookrightarrow\calV(P)$, i.e., $-T\beta(\id_M,p)(\xi,0)=\xi(p)$.
\end{corollary}
\begin{proposition}
	If all infinitesimal automorphisms in $\Kill(P,\kappa)$ are complete, then the group $H:=\langle\flow^{\xi}_1\colon \xi\in\Kill(P,\kappa)\rangle_{\Diff(P)}$ is a normal subgroup of $\Aut(P,\kappa)$. There is then a unique Banach--Lie group structure on $\Aut(P,\kappa)$ that makes $H$ an open Lie subgroup. The natural map $\tau\colon\Aut(P,\kappa)\times P\rightarrow P$ is a smooth action whose derived action is the inclusion map $\Kill(P,\kappa)\hookrightarrow\calV(P)$.
\end{proposition}
\begin{proof}
	The group $H$ is a subgroup of $\Aut(P,\kappa)$ by definition. We shall show that it is normal. Given any $g\in \Aut(P,\kappa)$, we have to check that $gHg^{-1}\subseteq H$. The set of all $h\in H$ satisfying $ghg^{-1}\subseteq H$ is a subgroup of $H$. Therefore, it suffices to verify $g\flow^{\xi}_1g^{-1}\subseteq H$ for all $\xi\in\Kill(P,\kappa)$. By Proposition~\ref{prop:fRelatedVectorFields}, we have $g\flow^{\xi}_tg^{-1}=\flow^{g_\ast\xi}_t$ for all $t\in\RR$. As then the flow maps $\flow^{g_\ast\xi}_t$ are in  $\Aut(P,\kappa)$, the vector field $g_\ast\xi$ is in $\Kill(P,\kappa)$, so that $g\flow^{\xi}_1g^{-1}\in H$.
	
	To show the existence of a Lie group structure, we have to check that for each $g\in \Aut(P,\kappa)$, the restriction $c_g|_H$ of the conjugation map $c_g$ is a smooth automorphism of $H$ (cf.\ \cite[Cor.~II.2.3]{Nee06}).
	As it is a homomorphism, it suffices to verify the smoothness in a neighborhood of $\id_P$. For this, we work in exponential charts.
	By the preceding considerations, we have $c_g|_H\circ\exp=\exp\circ g_\ast$, so that we shall show that the linear map $g_\ast\colon\Kill(P,\kappa)\rightarrow\Kill(P,\kappa)$ is continuous. Using an embedding $\ev_p\colon \Kill(P,\kappa)\hookrightarrow T_pP$, this follows from the continuity of $(T_pg)|_{(T_pP)^{H_p}}\colon (T_pP)^{H_p}\rightarrow (T_pP)^{H_p}$.
	
	The remaining statement follows from $\tau|_{H\times P}=\beta$ and $\tau|_{gH\times P}=g\circ\beta\circ(\lambda_{g^{-1}}|_{gH}\times \id_P)$ for all $g\in \Aut(P,\kappa)$, where $\lambda_{g^{-1}}$ denotes the left multiplication with $g^{-1}$ in $\Aut(P,\kappa)$.
\end{proof}
\begin{proposition}
	Assuming all infinitesimal automorphisms in $\Kill(P,\kappa)$ to be complete, we turn $\Aut(P,\kappa)$ into a Banach--Lie group. Then for each $p\in P$, the orbit map\linebreak $\tau_p\colon\Aut(P,\kappa)\rightarrow P,\ g\mapsto g(p)$ is an injective local topological embedding.
\end{proposition}
\begin{proof}
	We shall show the injectivity of $\tau_p$. Given any $g_1,g_2\in\Aut(P,\kappa)$ with $g_1(p)=g_2(p)$, the automorphism $g:=g_2^{-1}g_1$ satisfies $g(p)=p$. We have to check that $g=\id_P$ and shall do this by showing that the fixed point set $\Fix(g)$ is all of $P$. Given any $p^\prime\in P$,  there is a map $s\in S$ with $s(p)=p^\prime$ (cf.\ Lemma~\ref{lem:SconnectsM}).
	By Lemma~\ref{lem:fs=sf,sInS}, we have $g(s(p))=s(g(p))=s(p)$, so that $p^\prime\in\Fix(g)$, hence, $\Fix(g)=P$ follows.
	
	To see that $\tau_p$ is locally a topological embedding, it suffices to check this around $\id_P\in\Aut(P,\kappa)$, as
	$\tau_p= g\circ\tau_p\circ\lambda_{g^{-1}}$
	for all $g\in \Aut(P,\kappa)$.
	In view of Lemma~\ref{lem:T_xf=closedEmbedding}, it suffices to check that $T_{\id_P}\tau_p$ is a closed embedding. Indeed, $-T_{\id_P}\tau_p=\ev_p\colon \Kill(P,\kappa)\rightarrow T_pP$ is a closed embedding (cf.\ Corollary~\ref{cor:Banach-LieAlgebraKill(M,kappa,rho)}).
\end{proof}
%
%
%
%
%
%
\section{The Automorphism Group of an Affine Banach Manifold}
Given a connected affine Banach manifold $M$ that is geodesically complete, we show that its automorphism group can be turned into a Banach--Lie group acting smoothly on $M$.

In this section, we first collect a number of definitions and properties concerning affine connections on Banach manifolds. Given a connected affine manifold $(M,\nabla)$, the soldering form and the connection form equip the frame bundle $\Fr(M)$ with a $\{\1\}$-structure $\kappa$. A diffeomorphism $f$ of $M$ is affine if and only if its induced automorphism $\Fr(f)$ of the frame bundle leaves $\kappa$ invariant. We observe that the automorphism groups $\Aut(M,\nabla)$ and $\Aut(\Fr(M),\kappa)$ as well as the Lie algebras $\Kill(M,\nabla)$ and $\Kill(\Fr(M),\kappa)$ of infinitesimal automorphisms are naturally isomorphic, respectively. Assuming $(M,\nabla)$ to be geodesically complete, we show that these Lie algebras consist of complete vector fields. Therefore we can apply the results of the preceding section.

\subsection{Affine Connections on the Tangent Bundle} \label{sec:affineConnections}
Let $M$ be a smooth Banach manifold, $\pi\colon TM\rightarrow M$ be the natural projection of its tangent bundle and $\pi_{TM}\colon TTM \rightarrow TM$ be the natural projection of the tangent bundle of $TM$. Also the map $T\pi\colon TTM\rightarrow TM$ makes $TTM$ a vector bundle over $TM$ (cf.\ \cite[p.~104]{Lan01}). The composition $\pi\circ\pi_{TM}$ turns $TTM$ into a fiber bundle over $M$.

An \emph{affine connection on $TM$} is a morphism $B\colon TM\oplus TM \rightarrow TTM$ of fiber bundles over $M$ such that $(\pi_{TM},T\pi)\circ B=\id_{TM\oplus TM}$ and such that $B$ is bilinear, i.e., for each $x\in M$, $B_x\colon T_xM\oplus T_xM \rightarrow TTM$ is bilinear. Note that $B_x(v,\cdot)\colon T_xM \rightarrow T_v(TM)$ and $B_x(\cdot,w)\colon T_xM \rightarrow (T\pi)^{-1}(w)$ are indeed maps between Banach spaces. The pair $(M,B)$ is called an \emph{affine Banach manifold}.

In a chart $\varphi\colon U\rightarrow V\subseteq E$, an affine connection $B$ can be written as
$$\begin{array}{cccc}
	TT\varphi\circ B\circ (T\varphi\oplus T\varphi)^{-1}\colon & TV\oplus TV=V\times E\times E   & \rightarrow & TTV=V\times E\times E\times E\\
	& (x,v,w) & \mapsto & (x,v,w,B^\varphi_x(v,w))
\end{array}$$
with a smooth map
$B^\varphi\colon V \rightarrow \calL^2(E,E)$
from $V$ into the space of continuous bilinear maps $E\times E \rightarrow E$, which we call a \emph{local representation of $B$}. Considering two charts $\varphi_1$ and $\varphi_2$, the change of variable formula for the transition map $h:=\varphi_2\circ\varphi_1^{-1}$ is given by
$$B^{\varphi_2}_{h(x)}(dh(x)(v),dh(x)(w))=d^2h(x)(v,w)+dh(x)(B^{\varphi_1}_{x}(v,w)).$$

An affine connection can also be given by a \emph{covariant derivative} $\nabla$, i.e., by a collection $(\nabla^U)_{U \subseteq M \, \mbox{\scriptsize open}}$ of $\RR$-bilinear maps
$$\nabla^U\colon {\cal V}(U) \times {\cal V}(U) \rightarrow {\cal V}(U),\ (\xi,\eta) \mapsto (\nabla^U)_\xi\eta $$
satisfying the conditions
\begin{enumerate}
	\item[(1)] $(\nabla^U)_{f\xi}\eta = f(\nabla^U)_\xi\eta$ \quad ($C^\infty(U)$-linearity in the first variable)
	\item[(2)] $(\nabla^U)_\xi(f\eta) = (\xi.f)\eta + f(\nabla^U)_\xi\eta$ \quad (derivation property)
\end{enumerate}
for all $\xi,\eta\in{\cal V}(U)$ and smooth functions $f\in C^\infty(U)$ such that the maps $\nabla^U$ are compatible in the sense that
$$((\nabla^{U_1})_\xi\eta)|_{U_2} = (\nabla^{U_2})_{\xi|_{U_2}}\eta|_{U_2}$$
for all $\xi,\eta\in{\cal V}(U_1)$, $U_2\subseteq U_1\subseteq M$.
In the following, we shall often suppress the index set $U$ by writing $$\nabla\!_\xi\eta:= (\nabla^U)_\xi\eta$$
for all $\xi,\eta\in {\cal V}(U)$.

There is a one-to-one correspondence between affine connections and covariant derivatives. It is determined by the local formula
$$(\nabla\!_\xi\eta)^\varphi(x) = d\eta^\varphi(x)(\xi^\varphi(x)) - B^\varphi_x(\eta^\varphi(x),\xi^\varphi(x)),$$
where $(\nabla\!_\xi\eta)^\varphi$, $\eta^\varphi$ and $\xi^\varphi$ denote the local representations of the vector fields.
As far as the vector field $\xi$ is concerned, $(\nabla\!_\xi\eta)(x)$ only depends on $\xi(x)$. Therefore, it make sense to define $\nabla\!_v\eta$ for vectors $v$.

Given an affine connection, there is a unique vector bundle morphism $K\colon TTM \rightarrow TM$ (over $\pi$) between the vector bundles $\pi_{TM}\colon TTM\rightarrow TM$ and $\pi\colon TM\rightarrow M$, such that
$$\nabla\!_\xi\eta \ =\ K\circ T\eta\circ\xi$$
for all vector fields $\xi,\eta\in\calV(U)$ with an open submanifold $U$ of $M$. It is called the \emph{connector}. In a chart $\varphi\colon U\rightarrow V\subseteq E$, it can be written as
$$\begin{array}{cccc}
	T\varphi\circ K\circ (TT\varphi)^{-1}\colon & TTV=V\times E\times E\times E  & \rightarrow & TV=V\times E\\
	& (x,v,w,z) & \mapsto & (x,z-B^\varphi_x(v,w)).
\end{array}$$
Note that for each $v\in T_xM$ (with $x\in M$), a vector $z\in T_v(TM)$ can be given by $T\pi(z)$ and $K(z)$ and it is called \emph{vertical} if $T\pi(z)=0_{x}$ and \emph{horizontal} if $K(z)=0_x$.

Given a smooth curve $\alpha\colon J \rightarrow M$, let $\gamma\colon J \rightarrow TM$ be a lift of $\alpha$ to $TM$, i.e., a curve on $TM$ satisfying $\pi\circ\gamma=\alpha$. The {\em derivative of $\gamma$ along $\alpha$} is the unique lift $\nabla\!_{\alpha^\prime}\gamma$ of $\alpha$ to $TM$ that in a chart $\varphi\colon U \rightarrow V\subseteq E$ has the expression
$$(\nabla\!_{\alpha^\prime}\gamma)^\varphi(t) \ =\ (\gamma^\varphi)^\prime(t)-B^\varphi_{\alpha^\varphi(t)}(\gamma^\varphi(t),(\alpha^\varphi)^\prime(t)).$$
We also use the notation $\nabla\!_{\alpha^\prime(t)}\gamma$. A lift $\gamma$ of $\alpha$ is said to be {\em $\alpha$-parallel} if $\nabla\!_{\alpha^\prime}\gamma = 0$. Note that this means that all tangent vectors $\gamma^\prime(t)$ are horizontal.

An affine connection induces {\em parallel transport} along smooth curves. For a curve \linebreak $\alpha\colon J \rightarrow M$ and $t_0,t_1 \in J$, we denote it by
$$P_{t_0}^{t_1}(\alpha)\colon T_{\alpha(t_0)}M \rightarrow T_{\alpha(t_1)}M.$$
It is a topological linear isomorphism and is defined by the property that for each $v\in T_{\alpha(t_0)}M$, the map
$$\gamma_v:=P_{t_0}^{(\cdot)}(\alpha)(v)\colon J \rightarrow TM$$
is the unique curve in $TM$ that is $\alpha$-parallel and satisfies
$\gamma_v(t_0)=v$.
In any chart\linebreak $\varphi\colon U\rightarrow V\subseteq E$, it then satisfies the linear differential equation
$$(\gamma_v^\varphi)^\prime(t) \ =\ B^\varphi_{\alpha^\varphi(t)}(\gamma_v^\varphi(t),(\alpha^\varphi)^\prime(t))$$
and it is uniquely determined by satisfying these equations
for a collection of charts covering the curve $\alpha$ and by satisfying the initial condition $\gamma_v(t_0)=v$.
%

A {\em geodesic} is a curve $\alpha$ in $M$ whose derivative $\alpha^\prime$ is $\alpha$-parallel, i.e.,
$\nabla\!_{\alpha^\prime}\alpha^\prime = 0$.
For each $v\in T_xM,\, x\in M$, for which the unique maximal geodesic $\alpha_v\colon J \rightarrow TM$ with $\alpha_v^\prime(0)=v$ satisfies $1\in J$, we define
$$\exp(v):=\exp_x(v) := \alpha_v(1).$$
We denote the open domains of $\exp$ and $\exp_x$ by ${\cal D}_{\exp}\subseteq TM$ and ${\cal D}_{\exp_x}\subseteq T_xM$, respectively,
and get smooth maps $\exp\colon {\cal D}_{\exp} \rightarrow M$ and $\exp_x:=\exp|_{T_xM\cap \calD_{\exp}} \colon {\cal D}_{\exp_x} \rightarrow M$. Each geodesic $\alpha\colon J \rightarrow TM$ with $\alpha^\prime(0)=v$ satisfies $\alpha(t)=\exp(tv)$.
A manifold with an affine connection is called \emph{geodesically complete} if the domain of each maximal geodesic is all of $\RR$.

Let $V\subseteq {\cal D}_{\exp_x}$ be an open neighborhood of $0$ in $T_xM=:E$ that is star-shaped with respect to $0$ (i.e., $[0,1]V\subseteq V$) such that $\exp_x$ induces a diffeomorphism of $V$ onto its open\linebreak image $W$. Then $W$ is said to be a \emph{normal neighborhood of $x$}. We call the chart\linebreak $\varphi:=(\exp|_V^W)^{-1}\colon W \rightarrow V\subseteq E$ a \emph{normal chart at $x$}.
Normal neighborhoods do exist, as $\exp_x\colon {\cal D}_{\exp_x}\rightarrow M$ induces a local diffeomorphism at $0\in T_xM$, since $T_0\exp_x=\id_{T_xM}$ (cf.\ \cite[Th.~IV.4.1]{Lan01}).

Further details can be found in \cite[IV, VIII and X]{Lan01}, but basically for the case of torsionfree connections.
Cf.\ also \cite{KN63}, \cite{Kli82} and \cite{Ber08} for more material on connections.
\subsection{The Frame Bundle of an Affine Banach Manifold}
Let $M$ be a smooth Banach manifold (with model space $E$).
The set $$\Fr(M):=\bigcup_{x\in M}\Iso(E,T_xM)$$ (of topological linear isomorphisms) equipped with the projection $q\colon \Fr(M)\rightarrow M,\linebreak \Iso(E,T_xM)\ni p \mapsto x$ carries the structure of a smooth $\GL(E)$-principal bundle with respect to the action
$$\rho\colon\Fr(M) \times \GL(E)\rightarrow \Fr(M),\ (p,g)\mapsto p.g:=p\circ g.$$

More precisely, for each chart $\varphi\colon U\rightarrow V\subseteq E$ of $M$, the map
$$\begin{array}{ccccc}
	\Fr(\varphi)\colon & \Fr(U) &\rightarrow & V\times \GL(E)& \subseteq E\times\gl(E) \\
	&\Iso(E,T_xU)\ni p & \mapsto & (\varphi(x),d\varphi(x)\circ p).
\end{array}$$
is a bundle chart of $\Fr(M)$, and we have $$q(\Fr(\varphi)^{-1}(\varphi(x),g))=x \quad\mbox{and}\quad \Fr(\varphi)^{-1}(\varphi(x),g_1g_2)=\Fr(\varphi)^{-1}(\varphi(x),g_1).g_2$$
for all $x\in U$ and $g, g_1,g_2\in \GL(E)$. The bundle $\Fr(M)$ is called the \emph{frame bundle over $M$}. For further details, see \cite[7.10.1]{Bou07}.

The \emph{soldering form} $\theta$ on $\Fr(M)$ is the $E$-valued 1-form on $\Fr(M)$ defined by
$$\theta_p\colon T_p\Fr(M)\rightarrow E,\ v\mapsto p^{-1}(Tq(v)).$$
With respect to a chart $\varphi\colon U\rightarrow V\subset E$ of $M$ (and the corresponding bundle chart $\Fr(\varphi)$), its local representation $\theta^\varphi\colon V\times\GL(E)\rightarrow \calL(E\times\gl(E),E)$ is given by
$\theta^\varphi_{(x,g)}(v,w)  = g^{-1}(v)$
with $(x,g)\in V\times \GL(E)$ and $(v,w)\in E\times\gl(E)$.
 
Given an affine connection on $TM$, the \emph{connection form} $\omega$ on $\Fr(M)$ is the $\gl(E)$-valued 1-form on $\Fr(M)$ defined by
$$\omega_p\colon T_p\Fr(M)\rightarrow \gl(E),\ \omega_p(v)(e)= p^{-1}\big(K(T\widehat e(v))\big),$$ where for each $e\in E$, the map $\widehat e\colon \Fr(M)\rightarrow TM$ is the bundle morphism over $M$ given by $\widehat e(p)=p(e)$.
With respect to a chart $\varphi\colon U\rightarrow V\subset E$ of $M$, the map $T\widehat e$ can be written as
$$\begin{array}{cccc}
	TT\varphi \circ T\widehat e\circ T\Fr(\varphi)^{-1}\colon & V\times\GL(E)\times E\times\gl(E)&\rightarrow & V\times E \times E \times E \\
	& (x,g,v,w) & \mapsto & (x,g(e),v,w(e))
\end{array}$$
and the local representation $\omega^\varphi\colon V\times\GL(E)\rightarrow \calL(E\times\gl(E),\gl(E))$ of $\omega$ is given by
$$\omega^\varphi_{(x,g)}(v,w)(e) \ =\ g^{-1}(w(e)-B^\varphi_x(g(e),v))$$
with $(x,g)\in V\times\GL(E)$, $(v,w)\in E\times\gl(E)$ and $e\in E$.

For each $p\in\Fr(M)$, a vector $v\in T_p\Fr(M)$ can be given by $\theta_p(v)$ and $\omega_p(v)$ and it is called \emph{vertical} if $\theta_p(v)=0$ and \emph{horizontal} if $\omega_p(v)=0$. More precisely, the $(E\times\gl(E))$-valued 1-form $\kappa:=(\theta,\omega)$ is a $\{\1\}$-structure on $\Fr(M)$.
For each $e\in E$, the tangent map $T\widehat e\colon T\Fr(M)\rightarrow TTM$ maps vertical vectors in $T_p\Fr(M)$ to vertical ones in $T_{p(e)}TM$, and horizontal vectors to horizontal ones.

For each $\lambda\in E$, we define the \emph{standard horizontal vector field} $H_\lambda$ on $\Fr(M)$ by the requirement that $\theta(H_\lambda)(p)\equiv \lambda$ and $\omega(H_\lambda)(p)\equiv 0$, i.e., $H_\lambda(p):=\kappa_p^{-1}(\lambda,0)$. For a chart $\varphi\colon U\rightarrow V\subseteq E$ of $M$, its local representation $H^\varphi_\lambda\colon V\times \GL(E) \rightarrow E\times \gl(E)$ is given by
$$H_\lambda^{\varphi}(x,g) \ =\ \big(g(\lambda),\ e\mapsto B^\varphi_x(g(e),g(\lambda))\big).$$
By working in charts, we observe that $H_{g^{-1}(\lambda)}\circ \rho_g = T\rho_g\circ H_\lambda$ for all $\rho_g\colon\Fr(M)\rightarrow\Fr(M)$ with $g\in \GL(E)$, i.e., $H_\lambda$ and $H_{g^{-1}(\lambda)}$ are $\rho_g$-related.
\begin{proposition}\label{prop:standardHorizontal}
	Given $\lambda\in E$ and $g\in\GL(E)$, we have:
	\begin{enumerate}
		\item[\rm (1)]
		Every (maximal) integral curve of $H_\lambda$ is mapped under
		$\rho_g$ to a (maximal) integral curve of $H_{g^{-1}(\lambda)}$.
		\item[\rm (2)]
		Every (maximal) integral curve $\gamma$ of $H_\lambda$ is mapped under $q\colon \Fr(M)\rightarrow M$ to a (maximal) geodesic in $M$, and we have $(q\circ \gamma)^\prime(t)=\gamma(t)(\lambda)$.
		\item[\rm (3)]
		For each (maximal) geodesic in $M$, there is a  (maximal) integral curve of $H_\lambda$ that is mapped to it under $q\colon \Fr(M)\rightarrow M$.
	\end{enumerate}
	Given $p\in \Fr(M)$ with $q(p)=x$, we have:
	\begin{enumerate}
		\item[\rm (4)]
		For each (maximal) geodesic $\alpha$ in $M$ with initial condition $\alpha(0)=x$, there is a $\lambda\in E$ such that $H_\lambda$ possesses an (maximal) integral curve $\gamma$ with initial condition $\gamma(0)=p$ that is mapped to $\alpha$ by $q$.
	\end{enumerate}
\end{proposition}
\begin{proof}
	(1) is obvious by the fact that $H_\lambda$ and $H_{g^{-1}(\lambda)}$ are $\rho_g$-related.
	
	(2) Given a curve $\gamma\colon J\rightarrow \Fr(M)$, we represent it with respect to a given chart\linebreak $\varphi\colon U \rightarrow V\subseteq E$ of $M$ by the form $(\alpha^\varphi,\gamma^\varphi)\colon \alpha^{-1}(U)\rightarrow E\times \GL(E)$ with $\alpha:=q\circ \gamma$. If $\gamma$ is an integral curve of $H_\lambda$, then we have in a chart $$(\alpha^\varphi)^\prime(t)=\gamma^\varphi(t)(\lambda) \quad\mbox{and}\quad (\gamma^\varphi)^\prime(t)(e)=B^\varphi_{\alpha^\varphi(t)}(\gamma^\varphi(t)(e),\gamma^\varphi(t)(\lambda)).$$
	Therefore, we obtain $(\alpha^\varphi)^{\prime\prime}(t)=(\gamma^\varphi)^\prime(t)(\lambda)=B^\varphi_{\alpha^\varphi(t)}((\alpha^\varphi)^\prime(t),(\alpha^\varphi)^\prime(t))$, so that $\alpha$ is a geodesic with $\alpha^\prime(t)=\gamma(t)(\lambda)$.
	
	Assume now that $\gamma$ is maximal. Let $\beta\colon I\rightarrow M$ be the maximal geodesic with $J\subseteq I$ and $\beta|_J = \alpha$. To see that $J=I$, we assume (for contradiction) that there is a boundary point $t_1$ of $J$ in $I$.	
	With respect to a chart $\varphi\colon U\rightarrow V\subseteq E$ around $\beta(t_1)$, we have
	$$(\beta^\varphi)^\prime(t)=\gamma^\varphi(t)(\lambda) \quad\mbox{and}\quad (\gamma^\varphi)^\prime(t)=B^\varphi_{\beta^\varphi(t)}(\cdot,(\beta^\varphi)^\prime(t))\circ \gamma^\varphi(t)$$
	on $\alpha^{-1}(U)$. We shall extend $\gamma^\varphi$ around $t_1$ such that $(\beta^\varphi,\gamma^\varphi)$ still satisfies these differential equations, contradicting the maximality of $\gamma$. As the second equation is a linear differential equation for $\gamma^\varphi$, we can extend $\gamma^\varphi$ around $t_1$. To see that this extension also satisfies the first equation, we observe
	$$(\gamma^\varphi)^\prime(t)(\lambda) \ =\ B^\varphi_{\beta^\varphi(t)}(\gamma^\varphi(t)(\lambda),(\beta^\varphi)^\prime(t))$$
	and remember that this differential equation is also satisfied by $(\beta^\varphi)^\prime$, the curve $\beta$ being a geodesic. Hence, we have 
	$(\beta^\varphi)^\prime(t)=\gamma^\varphi(t)(\lambda)$ by uniqueness.
	
	(3) Given a geodesic $\beta\colon I\rightarrow M$, choose a frame $p_0\in\Fr(M)$ such that $p_0(\lambda)=\beta^\prime(t_0)$ for some $t_0\in I$. Let $\gamma\colon J\rightarrow \Fr(M)$ be the maximal integral curve of $H_\lambda$ with initial condition $\gamma(t_0)=p_0$. We know from (2) that the curve $\alpha:=q \circ \gamma$ is a maximal geodesic with $\alpha^\prime(t_0)=\gamma(t_0)(\lambda)=p_0(\lambda)=\beta^\prime(t_0)$, so that $I\subseteq J$ and $\alpha|_I=\beta$.	 
	
	(4) Choose some $\lambda_0\in E$. By (3),  we have a (maximal) integral curve $\gamma_0$ of $H_{\lambda_0}$ that is mapped to $\alpha$ under $q$. Let $g\in \GL(E)$ such that $\gamma_0(0).g=p$. Then the assertion follows by (1) when considering $\gamma:=\rho_g\circ\gamma_0$ and $\lambda:=g^{-1}(\lambda_0)$.
\end{proof}
\begin{corollary}\label{cor:completeH_lambda}
	Given $\lambda\in E$, the affine manifold $(M,\nabla)$ is geodesically complete if and only if the vector field $H_\lambda$ is complete.
\end{corollary}
\subsection{Affine Maps}\label{sec:affineMaps}
Given two affine manifolds $(M_1,B_1)$ and $(M_2,B_2)$, a map $f\colon M_1\rightarrow M_2$ is called \emph{affine}, if $TTf\circ B_1=B_2\circ (Tf\oplus Tf)$. Working with charts $\varphi_1\colon U_1\rightarrow V_1\subseteq E_1$ of $M_1$ and\linebreak $\varphi_2\colon U_2\rightarrow V_2\subseteq E_2$ of $M_2$ such that $f(U_1)=U_2$, this can be written as
\begin{equation}\label{eqn:affineMaps}
	d^2f^\varphi(x)(v,w)+df^\varphi(x)((B_1^{\varphi_1})_x(v,w))\ =\ (B_2^{\varphi_2})_{f^\varphi(x)}(df^\varphi(x)(v),df^\varphi(x)(w))
\end{equation}
for all $x$ in the domain of the local representation $f^\varphi\colon V_1\rightarrow V_2$ of $f$ and $v,w\in E_1$.
The following lemma is an easy result when taking a closer look at the corresponding local formulas.
\begin{lemma}
	Given a smooth map $f\colon M_1\rightarrow M_2$, the following are equivalent:
	\begin{enumerate}
		\item[\rm (a)]
		$f$ is affine.
		\item[\rm (b)]
		The connectors $K_1$ and $K_2$ are $Tf$-related in the sense that $Tf\circ K_1=K_2\circ TTf$.
		\item[\rm (c)]
		$TTf$ maps horizontal vectors to horizontal ones, i.e., $K_1(v)=0$ implies \linebreak $K_2(TTf(v))=0$ for all $v\in TTM$.
	\end{enumerate}
	(Note that vertical vectors are mapped to vertical ones regardless of whether $f$ is affine.)
\end{lemma}

Affine maps are compatible with parallel transport along curves, i.e., $T_{\alpha(t_1)}f\circ P_{t_0}^{t_1}(\alpha) = P_{t_0}^{t_1}(f\circ\alpha)\circ T_{\alpha(t_0)}f$ for all curves $\alpha\colon J\rightarrow M_1$ with $t_0,t_1\in J$. Geodesics  are mapped to geodesics. Further, we have $Tf(\calD_{\exp,1})\subseteq \calD_{\exp,2}$ and $f\circ\exp=\exp\circ Tf$.
A consequence is that, given an affine map, its values on connected components are uniquely determined by the tangent map at a single point, i.e., given affine maps $f,g\colon M_1\rightarrow M_2$ with $T_xf=T_xg$ for some $x\in M_1$, we have $f=g$ if $M_1$ is connected (cf.\ proof of \cite[Lem.~3.5]{Nee02}).

Affine maps are compatible with covariant derivatives of related vector fields, i.e., \linebreak $Tf(\nabla\!_v\eta_1)=\nabla\!_{Tf(v)}\eta_2$ for all $v\in TM_1$ and $\eta_1\in\calV(M_1)$, $\eta_2\in\calV(M_2)$ with $Tf\circ\eta_1=\eta_2\circ f$.

We now assume $M_1$ and $M_2$ to be modelled on the same Banach space $E$. A diffeomorphism $f\colon M_1\rightarrow M_2$ induces a principal bundle isomorphism $\Fr(f)\colon \Fr(M_1)\rightarrow \Fr(M_2)$ over $f$ defined by $\Fr(f)(p)=T_xf\circ p$ where $p\in \Iso(E,T_xM_1)$. It relates the soldering forms $\theta_1$ and $\theta_2$, i.e., $\Fr(f)^\ast\theta_2 =\theta_1$.

Conversely, every fiber-preserving diffeomorphism $F\colon\Fr(M_1)\rightarrow \Fr(M_2)$ with $F^\ast\theta_2 =\theta_1$ is induced by a unique diffeomorphism $f\colon M_1\rightarrow M_2$. Indeed, being fiber-preserving, $F$ induces a diffeomorphism $f$ between $M_1$ and $M_2$ and, by means of $F^\ast\theta_2 =\theta_1$, we can deduce\footnote
{From $(\theta_1)_p(v)=(\theta_2)_{F(p)}(TF(v))$ for all $p\in\Fr(M_1)$ and $v\in T_p\Fr(M_1)$, we obtain $p^{-1}(Tq(v))=\linebreak (Tf^{-1}\circ F(p))^{-1}(Tq(v))$, which leads to $Tf\circ p=F(p)$, i.e., $\Fr(f)(p)=F(p)$.}
$\Fr(f)=F$.

\begin{lemma}\label{lem:affineMapsRelateConnectionForms}
	Given a diffeomorphism $f\colon M_1\rightarrow M_2$ between affine manifolds, the following are equivalent:
	\begin{enumerate}
		\item[\rm (a)]
		$f$ is affine.
		\item[\rm (b)]
		The connection forms $\omega_1$ and $\omega_2$ are $\Fr(f)$-related, i.e., $\Fr(f)^\ast\omega_2=\omega_1$.
		\item[\rm (c)]
		$T\Fr(f)$ maps horizontal vectors to horizontal ones, i.e.,
		$$(\omega_1)_p(v)=0 \ \Rightarrow\ (\omega_2)_{\Fr(f)(p)}(T_p\Fr(f)(v))=0$$
		for all $p\in \Fr(M_1)$ and $v\in T_p\Fr(M_1)$.
	\end{enumerate}
	(Note that vertical vectors are mapped to vertical ones regardless of whether $f$ is affine.)
\end{lemma}
\begin{proof}
	We first note that $Tf\circ \widehat e^{\scriptscriptstyle 1}=\widehat e^{\scriptscriptstyle 2}\circ\Fr(f)$ for all $e\in E$.
	
	(a)$\Rightarrow$(b): As $Tf$ relates the connectors $K_1$ and $K_2$, we have
	\begin{eqnarray*}
		(\omega_2)_{\Fr(f)(p)}(T_p\Fr(f)(v))(e) &=& (\Fr(f)(p))^{-1}\big(K_2(T\widehat e^{\scriptscriptstyle 2}(T_p\Fr(f)(v)))\big) \\
		&=& (T_xf\circ p)^{-1}\big(T_xf(K_1(T\widehat e^{\scriptscriptstyle 1}(v)))\big) = p^{-1}\big(K_1(T\widehat e^{\scriptscriptstyle 1}(v))\big) \\
		&=& (\omega_1)_p(v)(e)
	\end{eqnarray*}
	for all $p\in\Iso(E,T_xM_1)\subseteq\Fr(M_1)$, $v\in T_p\Fr(M_1)$ and $e\in E$.
	
	(b)$\Rightarrow$(c) is obvious.
	
	(c)$\Rightarrow$(a): Working with charts, we shall show (\ref{eqn:affineMaps}). We choose some $g\in\GL(E)$ and $e\in E$ such that $v=g(e)$. We put $z:=(B^{\varphi_1}_1)_x(\cdot,w)\circ g\in\gl(E)$ and observe that then $(\omega^{\varphi_1}_1)_{(x,g)}(w,z)=0$. By assumption (c), we then know that
	$$(\omega^{\varphi_2}_2)_{(f^\varphi(x),df^\varphi(x)\circ g)}(df^\varphi(x)(w), d^2f^\varphi(x)(\cdot,w)\circ g + df^\varphi(x)\circ z) \ =\ 0,$$
	i.e.,
	$$d^2f^\varphi(x)(\cdot,w)\circ g + df^\varphi(x)\circ z \ =\ (B^{\varphi_2}_2)_{f^\varphi(x)}(\cdot,df^\varphi(x)(w))\circ (df^\varphi(x)\circ g).$$
	Applying both sides of this equation to $e$, we get (\ref{eqn:affineMaps}).
\end{proof}
\begin{proposition}\label{prop:isomorphismOfAut}
	Given an affine manifold $(M,\nabla)$, the map $\Aut(M,\nabla)\rightarrow \Aut(\Fr(M),\kappa),\linebreak f\mapsto \Fr(f)$ with $\kappa=(\theta,\omega)$ is an isomorphism of groups.
\end{proposition}
\begin{proof}
	The map is correctly defined, as for each $f\in \Aut(M,\nabla)$, the induced map $\Fr(f)$ is a principal bundle automorphism that preserves $\theta$ and $\omega$. It is bijective, since each $F\in\Aut(\Fr(M),\kappa)$ is induced by a unique diffeomorphism $f$ of $M$, which is affine by Lemma~\ref{lem:affineMapsRelateConnectionForms}. The map is a group homomorphism, as
	$$\Fr(f\circ g)(p) \ =\ T(f\circ g)\circ p\ =\ Tf\circ\Fr(g)(p) \ =\ \Fr(f)(\Fr(g)(p)) \ =\ (\Fr(f)\circ\Fr(g))(p)$$
	for all $f,g\in\Aut(M,\nabla)$.
\end{proof}
\subsection{Infinitesimal Affine Automorphisms}
Let $M$ be a smooth Banach manifold (with model space $E$)
and $q\colon \Fr(M)\rightarrow M$ the natural projection of the frame bundle.

Given $\xi\in\calV(M)$, the \emph{natural lift} of $\xi$ is defined as the vector field $\overline\xi\in\calV(\Fr(M))$ given by
$\overline\xi(p):=\left.\frac{d}{dt}\right|_{t=0}\Fr(\flow^{\xi}_t)(p)$.
The results of this section can essentially be found in \cite[VI.2]{KN63}  for the finite-dimensional case.
\begin{lemma}\label{lem:naturalLiftOnFrameBundle}
	Given $\xi\in\calV(M)$, we have:
	\begin{enumerate}
		\item[\rm (1)]
		$\overline\xi$ and $\xi$ are $q$-related, i.e., $Tq\circ \overline\xi= \xi\circ q$.
		\item[\rm (2)] Given a chart $\varphi\colon U\rightarrow V\subseteq E$ of $M$, the local representation $\overline\xi^{\varphi}\colon V\times \GL(E)\rightarrow E\times \gl(E)$ of $\overline\xi$ is given by
		$$\overline\xi^{\varphi}(x,g) \ =\ (\xi^\varphi(x),d\xi^\varphi(x)\circ g).$$
		\item[\rm (3)]
		The domain of the local flow $\flow^{\overline\xi}$ of $\overline\xi$ is given by
		$${\cal D}(\overline\xi) \ =\ \{(t,p)\in \RR\times \Fr(M)\colon (t,q(p))\in {\cal D}(\xi)\}.$$
		Each flow map $\flow^{\overline\xi}_t$ is the principal bundle isomorphism between the frame bundles $\Fr(\calD_t(\xi))$ and $\Fr(\calD_{-t}(\xi))$ induced by the flow map $\flow^\xi_t$. We thus have
		$$\flow^{\overline\xi}_t(p)=\Fr(\flow^\xi_t)(p) \quad\mbox{and}\quad
			\flow^\xi_t(q(p))= q(\flow^{\overline\xi}_t(p))$$
		for all $(t,p)\in {\cal D}(\xi)$.
		The second equation says that $q$ maps the maximal integral curve of $\overline\xi$ with initial value $p$ to the maximal integral curve of $\xi$ with initial value $q(p)$.
		\item[\rm (4)]
		$\overline\xi$ is invariant under all $\rho_g\colon \Fr(M)\rightarrow\Fr(M)$ with $g\in\GL(E)$, i.e., $(\rho_g)_\ast\overline\xi=\overline\xi$.
		\item[\rm (5)]
		$\lied_{\overline\xi}\theta=0$.
	\end{enumerate}
\end{lemma}
\begin{proof}
	(1) For all $p\in \Fr(M)$, we have
	\begin{eqnarray*}
		(Tq\circ\overline\xi)(p)
		&=& \textstyle Tq(\left.\frac{d}{dt}\right|_{t=0}\Fr(\flow^\xi_t)(p))
		\ = \ \textstyle\left.\frac{d}{dt}\right|_{t=0}q(\Fr(\flow^\xi_t)(p)) \\
		&=& \textstyle\left.\frac{d}{dt}\right|_{t=0}\flow^\xi_t(q(p))
		\ = \ \xi(q(p)).
	\end{eqnarray*}
	
	(2) Given any $(x,g)\in V\times \GL(E)$, let $\varphi_1\colon U_1\rightarrow V_1\subseteq E$ be a restriction of $\varphi$ and $\varepsilon>0$ such that $x\in V_1$, $W:={]-\varepsilon,\varepsilon[}\times U_1\subseteq {\cal D}(\xi)$ and $\flow^\xi(W)\subseteq U$.
	We deduce from the definition of the natural lift that the local representation $\overline\xi^{\varphi}$ of $\overline\xi$ is given by
	\begin{equation*}
		\textstyle \overline\xi^{\varphi}(x,g)
			= (\left.\frac{d}{dt}\right|_{t=0}(\flow^\xi_t)^{\varphi_1}(x),
				\left.\frac{d}{dt}\right|_{t=0}d(\flow^\xi_t)^{\varphi_1}(x)\circ g),
	\end{equation*}
	where
	$$(\flow^\xi_t)^{\varphi_1}:=(\flow^\xi_t)^{\varphi_1,\varphi}
		\colon V_1\rightarrow V$$
	denotes the local representation of $\flow^\xi_t$. 
	As $\left.\frac{d}{dt}\right|_{t=0}(\flow^\xi_t)^{\varphi_1}(x) = \xi^\varphi(x)$, it remains to verify that
	$$\left.\frac{d}{dt}\right|_{t=0}d(\flow^\xi_t)^{\varphi_1}(x)\circ g
		\ =\ d\xi^\varphi(x)\circ g.$$
	Writing $d(\flow^\xi_t)^{\varphi_1}(x)=d_2(\flow^\xi|_W)^{\varphi_1}(t,x)$,
	where 
	$$(\flow^\xi|_W)^{\varphi_1}
		:=(\flow^\xi|_W)^{(\id_{]-\varepsilon,\varepsilon[}\times\varphi_1),\varphi}\colon
			]-\varepsilon,\varepsilon[ \times V_1 \rightarrow V$$
	denotes the local representation of $\flow^\xi|_W$, we can use Schwarz's theorem, which carries over to Banach spaces (cf.\ \cite[Th.~XIII.7.3]{Lan93}), i.e., partial derivatives commute. Hence, the equation easily follows.
	
	(3) For each $p\in \Fr(M)$, the maximal integral curve $\flow^{\overline\xi}_p:=\flow^{\overline\xi}(\cdot,p)\colon J_p\rightarrow \Fr(M)$ is mapped by $q$ to an integral curve of $\xi$ with initial value $q(p)$, as
	\begin{equation*}
		(q\circ\flow^{\overline\xi}_p)^\prime(t)
			\ =\ Tq((\flow^{\overline\xi}_p)^\prime(t))
			\ =\ Tq\big(\overline\xi(\flow^{\overline\xi}_p(t))\big)
			\ =\ \xi\big(q(\flow^{\overline\xi}_p(t))\big).
	\end{equation*}
	Hence, $q\circ\flow^{\overline\xi}_p$ is a restriction of the maximal integral curve $\flow^\xi_{q(p)}\colon J_{q(p)}\rightarrow M$. To see that $q\circ\flow^{\overline\xi}_p=\flow^\xi_{q(p)}$, i.e., $J_p=J_{q(p)}$, we shall show that $J_{q(p)}\ni t \mapsto \Fr(\flow^\xi_t)(p)$ is an integral curve of $\overline\xi$. Indeed, we have
	\begin{eqnarray*}
		\left.\frac{d}{dt}\right|_{t=t_0}\Fr(\flow^\xi_t)(p)
		&=& \left.\frac{d}{dt}\right|_{t=0}\Fr(\flow^\xi_{t_0+t})(p)
		\ =\ \left.\frac{d}{dt}\right|_{t=0}\Fr(\flow^\xi_t)(\Fr(\flow^\xi_{t_0})(p)) \\
		&=& \overline\xi(\Fr(\flow^\xi_{t_0})(p)).
	\end{eqnarray*}
	From the preceding considerations, the assertions immediately follow.
	
	(4) By Proposition~\ref{prop:fRelatedVectorFields}, it suffices to check $\rho_g(\calD_t(\overline\xi))\subseteq\calD_t(\overline\xi)$ and $\rho_g\circ\flow^{\overline\xi}_t=\flow^{\overline\xi}_t\circ\rho_g|_{\calD_t(\overline\xi)}$.
	This is evident, since we have $\calD_t(\overline\xi)=q^{-1}(\calD_t(\xi))=\Fr(\calD_t(\xi))$ and since, being a principal bundle isomorphism, $\flow^{\overline\xi}_t=\Fr(\flow^\xi_t)\colon \Fr(\calD_t(\xi))\rightarrow\Fr(\calD_{-t}(\xi))$ intertwines the natural $\GL(E)$-actions on $\Fr(\calD_t(\xi))$ and $\Fr(\calD_{-t}(\xi))$, which are restrictions of the action $\rho$ on $\Fr(M)$.
	
	(5) By Proposition~\ref{prop:lieXiOmega=0}, it suffices to check that $\theta$ is invariant under all flow maps $\flow^{\overline\xi}_t$, i.e., $(\flow^{\overline\xi}_t)^\ast(\theta|_{\calD_{-t}(\overline\xi)})=\theta|_{\calD_{t}(\overline\xi)}$. This is evident, since $\flow^{\overline\xi}_t=\Fr(\flow^{\xi}_t)$ relates the soldering forms on $\Fr(\calD_t(\xi))$ and $\Fr(\calD_{-t}(\xi))$, which are restrictions of the soldering form $\theta$ on $\Fr(M)$.
\end{proof}
\begin{lemma}\label{lem:vectorFieldsBeingNaturalLifts}
	Each vector field $X\in\calV(\Fr(M))$ that satisfies $(\rho_g)_\ast X=X$ for all $g\in\GL(E)$ and $\lied_X\theta= 0$ is the natural lift of a unique vector field $\xi\in\calV(M)$, i.e., $X=\overline\xi$. The vector field $\xi$ is uniquely determined by the property that $X$ and $\xi$ are $q$-related.
\end{lemma}
\begin{proof}
	The proof is not more difficult than in the finite-dimensional case (cf.\ \cite[p.~229]{KN63}).
	Due to $(\rho_g)_\ast X=X$, the map $Tq\circ X$ factorizes over $q$, since
	$$Tq\circ X \circ \rho_g \ =\ Tq \circ T\rho_g \circ X \ =\ T(q \circ\rho_g) \circ X \ =\ Tq \circ X$$ for all $g\in\GL(E)$. Hence, there exists a vector field $\xi\in\calV(M)$ with $\xi\circ q=Tq\circ X$. 
	By Proposition~\ref{prop:fRelatedVectorFields}, we have
	$q(\calD_t(X))\subseteq\calD_t(\xi)$
	and
	$q\circ\flow^X_t=\flow^\xi_t\circ q|_{\calD_t(X)}$. To see that each flow map $\flow^X_t$ is equal to $\Fr(\flow^\xi_t)$, we have to check that $\theta$ is invariant under $\flow^X_t$ (cf.\ Section~\ref{sec:affineMaps}), but this is satisfied, due to Proposition~\ref{prop:lieXiOmega=0}. From $\flow^X_t=\Fr(\flow^\xi_t)=\flow^{\overline\xi}_t$ for all $t\in\RR$, we immediately obtain $X=\overline\xi$. The uniqueness is clear by Lemma~\ref{lem:naturalLiftOnFrameBundle}(1), as $q$ is surjective.
\end{proof}
From now on, let $(M,\nabla)$ ($=(M,B)$) be an affine manifold.
\begin{definition}
	A vector field $\xi\in\calV(M)$ is called an \emph{infinitesimal affine automorphism} if each flow map $\flow^\xi_t$ is an affine isomorphism. We denote the set of all infinitesimal automorphisms by $\Kill(M,\nabla)$ or $\Kill(M,B)$.
\end{definition}
\begin{lemma}\label{lem:infinitesimalAffineAutomorphisms}
	Given a vector field $\xi\in\calV(M)$, the following are equivalent:
	\begin{enumerate}
		\item[\rm (a)]
		$\xi$ is an infinitesimal affine automorphism.
		\item[\rm (b)]
		$\lied_{\overline\xi}\omega=0$.
		\item[\rm (c)]
		$\overline\xi\in\Kill(\Fr(M),\kappa)$, where $\kappa=(\theta,\omega)$.
		\item[\rm (d)]
		$[\overline\xi,H_\lambda]=0$ in the Lie algebra $\calV(\Fr(M))$ for all $\lambda\in E$.
	\end{enumerate}
\end{lemma}
\begin{proof}
	Cf.\ \cite[p.~230]{KN63} for the finite-dimensional case.
	
	(a)$\Rightarrow$(b): Every flow map $\flow^\xi_t$ is an affine isomorphism, so that the connection form $\omega$ is invariant by all induced maps $\Fr(\flow^\xi_t)=\flow^{\overline\xi}_t$ (cf.\ Lemma~\ref{lem:affineMapsRelateConnectionForms} and Lemma~\ref{lem:naturalLiftOnFrameBundle}). When applying Proposition~\ref{prop:lieXiOmega=0}, we get (b).
	
	(b)$\Rightarrow$(c): Together with $\lied_{\overline\xi}\theta=0$ (cf.\ Lemma~\ref{lem:naturalLiftOnFrameBundle}), we obtain $\lied_{\overline\xi}\kappa=0$ for the $\{\1\}$-structure $\kappa=(\theta,\omega)$. Further, we have $(\rho_g)_\ast\overline\xi=\overline\xi$ for all $g\in\GL(E)$, so that $\overline\xi\in\Kill(\Fr(M),\kappa)$ (cf.\ Section~\ref{sec:subsectionIntro}).
	
	(c)$\Rightarrow$(d):
	With respect to the notion of Section~\ref{sec:infinitesimalAutomorphisms}, we have $H_\lambda=\eta_{(\lambda,0)}$. Hence, (d) follows by Lemma~\ref{lem:Kill(M,kappa,rho)}.
	
	(d)$\Rightarrow$(a): By Lemma~\ref{lem:affineMapsRelateConnectionForms}(c), it suffices to check that for each flow map $\flow^\xi_t$, the map $T\Fr(\flow^\xi_t)=T\flow^{\overline\xi}_t$ maps horizontal vectors to horizontal ones. Given any $p\in\Fr(M)$ and $v\in T_p\Fr(M)$ with $\omega_p(v)=0$, we put $\lambda:=\theta_p(v)$. By the definition of $H_\lambda$, we then have $v=H_\lambda(p)$. By Corollary~\ref{cor:commutingFlows}, $H_\lambda$ is invariant under all flow maps $\flow^{\overline\xi}_t$, so that
	$$\omega_{\flow^{\overline\xi}_t(p)}(T\flow^{\overline\xi}_t(v)) \ =\ \omega_{\flow^{\overline\xi}_t(p)}(T\flow^{\overline\xi}_t(H_\lambda(p))) \ =\ \omega_{\flow^{\overline\xi}_t(p)}(H_\lambda(\flow^{\overline\xi}_t(p)))=0,$$
	as $\omega(H_\lambda)$ is identically 0.
\end{proof}
\begin{remark}\label{rem:infAffAutoChart}
	We can express Condition (d) in a chart $\varphi\colon U\rightarrow V\subseteq E$ of $M$ by using the local representations $\overline\xi^\varphi$ and $H_\lambda^\varphi$. A simple computation shows that
	$$dH_\lambda^\varphi(x,g)(\overline\xi^\varphi(x,g))\ = \ d\overline\xi^\varphi(x,g)(H_\lambda^\varphi(x,g)) \quad\mbox{for all $(x,g)\in V\times\gl(E)$ and $\lambda\in E$}$$
	is equivalent to
	\begin{eqnarray*}
		\lefteqn{d^2\xi^\varphi(x)(v,w)+d\xi^\varphi(x)(B^\varphi_x(v,w))} \\
		&=& dB^\varphi(x)(\xi^\varphi(x))(v,w)+B^\varphi_x(d\xi^\varphi(x)(v),w) + B^\varphi_x(v,d\xi^\varphi(x)(w))
	\end{eqnarray*}
	for all $x\in V$ and $v,w\in E$.
\end{remark}
\begin{proposition}\label{prop:isomorphismOfKill}
	The set $\Kill(M,\nabla)$ is a Lie subalgebra of $\calV(M)$ and the map
	$$\Kill(M,\nabla)\rightarrow \Kill(\Fr(M),\kappa)\subseteq \calV(\Fr(M)),\ \xi\mapsto \overline\xi$$
	is an isomorphism of Lie algebras. Hence, if $M$ is connected, $\Kill(M,\nabla)$ inherits the structure of a Banach--Lie algebra via this isomorphism (cf.\ Proposition~\ref{prop:Banach-LieAlgebraKill(M,kappa,rho)}).
	Its Banach space structure is then uniquely determined by the requirement that for each $p\in\Fr(M)$, the map
	$$\Kill(M,\nabla)\rightarrow T_p(\Fr(M)),\ \xi\mapsto \overline\xi(p)=\left.\frac{d}{dt}\right|_{t=0}\Fr(\flow^{\xi}_t)(p)$$
	is a closed embedding.
\end{proposition}
\begin{proof}
	Note that the map is correctly defined by Lemma~\ref{lem:infinitesimalAffineAutomorphisms} and that $\Kill(\Fr(M),\kappa)$ is a Lie algebra by Lemma~\ref{lem:Kill(M,kappa,rho)}. The map is bijective by Lemma~\ref{lem:vectorFieldsBeingNaturalLifts}. It suffices to show that $\Kill(\Fr(M),\kappa)\rightarrow \calV(M),\ \overline\xi\mapsto \xi$ is a homomorphism of Lie algebras.
	
	Let $\overline\xi_1$ and $\overline\xi_2$ be in $\Kill(\Fr(M),\kappa)$ and $\lambda\in\RR$. As $q_\ast\overline\xi_i=\xi_i$ for $i=1,2$, i.e., $\overline\xi_i$ and $\xi_i$ are $q$-related, we also have $q_\ast(\lambda\overline\xi_1+\overline\xi_2)=\lambda\xi_1+\xi_2$
	and $q_\ast[\overline\xi_1,\overline\xi_2]=[\xi_1,\xi_2]$ by the naturality of the Lie bracket. From Lemma~\ref{lem:vectorFieldsBeingNaturalLifts}, we then know that $\overline{\lambda\xi_1+\xi_2}=\lambda\overline\xi_1+\overline\xi_2$ and $\overline{[\xi_1,\xi_2]}=[\overline\xi_1,\overline\xi_2]$.
\end{proof}
\begin{corollary}\label{cor:embeddingOfInfAffAuto}
	If $M$ is connected, then, for each $x\in M$, the map
	$$\Kill(M,\nabla)\rightarrow T_xM \times \gl(T_xM),\ \xi\mapsto \big(\xi(x),\ v\mapsto\nabla\!_v\xi\big)$$
	is a closed embedding of Banach spaces.
\end{corollary}
\begin{proof}
	We choose some frame $p\in\Iso(E,T_xM)\subseteq \Fr(M)$. It suffices to construct an isomorphism $\Phi\colon T_p(\Fr(M))\rightarrow T_xM\times\gl(T_xM)$ satisfying $\Phi(\overline\xi(p))=(\xi(x),\nabla\!_\cdot\xi)$ for all $\xi\in\Kill(M,\nabla)$.
	Let $\varphi\colon U\rightarrow V\subseteq E$ be a chart at $x$ and put $(\bar x,g):=\Fr(\varphi)(p) =\linebreak (\varphi(x), d\varphi(x)\circ p) \in V\times\GL(E)$.
	We define 
	$$\Phi:=(d\varphi(x)\times d\varphi(x)_\ast)^{-1}\circ\Psi\circ d\Fr(\varphi)(p)$$ 
	with the isomorphisms
	$$\Psi\colon E\times\gl(E)\rightarrow E\times\gl(E),\ (v,w)\mapsto \big(v,\ w\circ g^{-1}-B^\varphi_{\bar x}(v,\cdot)\big)$$
	and	$d\varphi(x)_\ast\colon\gl(T_xM)\rightarrow\gl(E)$ given by $d\varphi(x)_\ast(w):=d\varphi(x)\circ w\circ d\varphi(x)^{-1}$.
	Given $\xi\in\Kill(M,\nabla)$, we have
	$$d\Fr(\varphi)(p)(\overline\xi(p)) \ =\ \overline\xi^\varphi(\bar x, g) \ =\ (\xi^\varphi(\bar x),d\xi^\varphi(\bar x)\circ g)$$
	(cf.\ Lemma~\ref{lem:naturalLiftOnFrameBundle}), which is mapped by $\Psi$ to $\big(\xi^\varphi(\bar x),d\xi^\varphi(\bar x)-B^\varphi_{\bar x}(\xi^\varphi(\bar x),\cdot)\big)$. Hence, $\Phi(\overline\xi(p))=(\xi(x),\nabla\!_\cdot\xi)$.
\end{proof}
\begin{lemma} \label{lem:epsilonComplete}
	Let $\xi$ be a smooth vector field on a smooth Banach manifold $N$ with local flow $\flow^\xi\colon\calD(\xi)\rightarrow N$. If there is an $\varepsilon>0$ such that ${[-\varepsilon,\varepsilon]}\times N \subseteq {\cal D}(\xi)$, then we have ${\cal D}(\xi)=\RR\times N$, i.e., $\xi$ is complete.
\end{lemma}
\begin{proof}
	Given an $\varepsilon>0$ such that ${[-\varepsilon,\varepsilon]}\times N \subseteq {\cal D}(\xi)$, we assume the opposite ${\cal D}(\xi)\neq\RR\times N$ and show that this will lead to a contradiction. In doing so, there is a greatest natural number $n$ satisfying
	$\textstyle {[-n\varepsilon,n\varepsilon]}\times N \subseteq {\cal D}(\xi).$
	Consequently, for each flow line $\flow^\xi_x\colon J_x\rightarrow N$ with $x\in N$, we have
	${[-n\varepsilon,n\varepsilon]}\subseteq J_x$, but at least
	${[-\varepsilon,\varepsilon]}\subseteq J_x$.
	For each $x\in N$, we have
	$J_{\flow^\xi(\varepsilon,x)}=J_x-\varepsilon$, so that
	$$\textstyle J_x = J_{\flow^\xi(\varepsilon,x)} + \varepsilon
		\supseteq {[-n\varepsilon,n\varepsilon]}
		+ \varepsilon \supseteq {[0,(n+1)\varepsilon]}.$$
	Similarly, we have $J_x\supseteq {[-(n+1)\varepsilon,0]}$, hence
	$\textstyle J_x \supseteq
		 {[-(n+1)\varepsilon,(n+1)\varepsilon]}.$
	This contradicts the maximality of $n$.
\end{proof}
\begin{theorem}\label{th:Kill(M,nabla)Complete}
	If $(M,\nabla)$ is geodesically complete, then all vector fields in $\Kill(M,\nabla)$ and in $\Kill(\Fr(M),\kappa)$ are complete.
\end{theorem}
\begin{proof}
	The proof is not more demanding than in the finite-dimensional case (cf.\ \cite[p.~234]{KN63}). By Lemma~\ref{lem:naturalLiftOnFrameBundle}(3) and Proposition~\ref{prop:isomorphismOfKill}, it suffices to prove the completeness of the vector fields in $\Kill(M,\nabla)$.
	Without loss of generality, we assume $M$ to be connected, as the matter of local flows and geodesics takes place in connected components. Given an infinitesimal automorphism $\xi\in\Kill(M,\nabla)$, it suffices to check ${[-\varepsilon,\varepsilon]}\times M \subseteq {\cal D}(\xi)$ for an $\varepsilon>0$ by Lemma~\ref{lem:epsilonComplete}.
	
	We consider some point $x_0\in M$ and let $\varepsilon>0$ be such that ${[-\varepsilon,\varepsilon]}\times \{x_0\} \subseteq {\cal D}(\xi)$. It suffices to show that the set $A$ of all $x\in M$ with ${[-\varepsilon,\varepsilon]}\times \{x\} \subseteq {\cal D}(\xi)$ is all of $M$. As $M$ is connected and $A$ is not empty, it suffices to check that $A$ and its complement $A^c$ both are open.
	
	To see that $A$ is open, let $x$ be any point in $A$. Due to Lemma~\ref{lem:naturalLiftOnFrameBundle}(3), we then have ${[-\varepsilon,\varepsilon]}\times \Iso(E,T_xM) \subseteq {\cal D}(\overline\xi)$.
	We consider a normal neighborhood $W$ of $x$ and shall show that
	it lies in $A$ and that therefore $A$ is open. Given any $y\in W$, we have to check that ${[-\varepsilon,\varepsilon]}\times \{y\} \subseteq {\cal D}(\xi)$. By Lemma~\ref{lem:naturalLiftOnFrameBundle}(3), it suffices to check ${[-\varepsilon,\varepsilon]}\times \{r\} \subseteq {\cal D}(\overline\xi)$ for some frame $r\in\Fr(M)$ at $y=q(r)$.
	The set $W$ being a normal neighborhood of $x$, there is a geodesic in $W$ that joins $x$ with $y$. Hence it follows by Proposition~\ref{prop:standardHorizontal} that there is a standard horizontal vector field $H_\lambda$ with $\lambda\in E$ and an integral curve $\gamma$ of $H_\lambda$ that joins some frame $p:=\gamma(0)$ at $x$ with some frame $r:=\gamma(1)$ at $y$.
	By Lemma~\ref{lem:infinitesimalAffineAutomorphisms}, we have $[\overline\xi,H_\lambda]=0$, so that we can apply Corollary~\ref{cor:commutingFlows}. As $H_\lambda$ is complete (cf.\ Corollary~\ref{cor:completeH_lambda}) and as we have $r=\flow^{H_\lambda}_1(p)$ and ${[-\varepsilon,\varepsilon]}\times \{p\} \subseteq {\cal D}(\overline\xi)$, we obtain ${[-\varepsilon,\varepsilon]}\times \{r\} \subseteq {\cal D}(\overline\xi)$.
	
	To see that $A^c$ is open, let $x$ be any point in $A^c$. We consider a normal neighborhood $W$ of $x$ and shall show that $W\subseteq A^c$. For each $y\in W$, there is a geodesic joining $y$ and $x$. If $y$ was in $A$ then $x$ would be in $A$, too, by the above argument. That is why $y$ is in $A^c$, hence $W\subseteq A^c$.
\end{proof}
\subsection{The Automorphism Group of a Geodesically Complete Affine Manifold}
\begin{theorem}\label{th:autoGroupOfAffineManifold}
	Let $(M,\nabla)$ be a connected affine Banach manifold that is geodesically complete.
	The automorphism group $\Aut(M,\nabla)$ can be turned into a Banach--Lie group such that
	$$\exp\colon \Kill(M,\nabla)\rightarrow \Aut(M,\nabla),\ \xi \mapsto \flow^{-\xi}_1$$
	is its exponential map. The natural map $\tau\colon \Aut(M,\nabla)\times M \rightarrow M$ is a smooth action whose derived action is the inclusion map $\Kill(M,\nabla)\hookrightarrow\calV(M)$, i.e., $-T\tau(\id_M,x)(\xi,0)=\xi(x)$.
	For each $p\in \Fr(M)$, the map $\Aut(M,\nabla) \rightarrow \Fr(M),\ f\mapsto \Fr(f)(p)$ is an injective local topological embedding.
\end{theorem}
\begin{proof}
	The assertions follow by Theorem~\ref{th:lieGroupKill(M,kappa,rho)}, Proposition~\ref{prop:isomorphismOfAut} and Proposition~\ref{prop:isomorphismOfKill}. Indeed, the exponential map of $\Aut(\Fr(M),\kappa)$ maps $\overline\xi\in\Kill(\Fr(M),\kappa)$ to $\flow^{-\overline\xi}_1=\Fr(\flow^{-\xi}_1)$, so that $\Kill(M,\nabla)\ni\xi\mapsto \flow^{-\xi}_1$ is an exponential map of $\Aut(M,\nabla)$.
	
	The natural smooth action $\overline\tau\colon\Aut(\Fr(M),\kappa)\times \Fr(M)\rightarrow \Fr(M)$ induces the smooth action $\tau\colon \Aut(M,\nabla)\times M\rightarrow M$, since $f(q(p))=q(\Fr(f)(p))$ for all $f\in\Aut(M,\nabla)$ and $p\in\Fr(M)$. Given $p\in \Fr(M)$, for the orbit maps $\overline\tau_p\colon \Aut(\Fr(M),\kappa)\rightarrow \Fr(M)$
	and
	$\tau_{q(p)}\colon\Aut(M,\nabla)\rightarrow M$,
	we have 
	$\tau_{q(p)}(f)=q(\overline\tau_p(\Fr(f))$
	and thus
	$$-T_{\id_M}\tau_{q(p)}(\xi) \ =\ -Tq(T_{\id_{\Fr(M)}}\overline\tau_p(\overline\xi)) \ =\ -Tq(-\overline\xi(p)) \ =\ \xi(q(p)).$$
	
	As the map $\Aut(M,\nabla)\rightarrow \Aut(\Fr(M),\kappa),\ f\rightarrow \Fr(f)$ is an isomorphism of Lie groups by construction, the map $\Aut(M,\nabla) \rightarrow \Fr(M),\ f\mapsto \Fr(f)(p)$ is an injective local topological embedding by Theorem~\ref{th:lieGroupKill(M,kappa,rho)}.
\end{proof}
\begin{acknowledgements}
	I am grateful to Karl-Hermann Neeb for his helpful communications and proof reading during my research towards this article. This work was supported by the Technical University of Darmstadt and by the Studienstiftung des deutschen Volkes.
\end{acknowledgements}
\bibliography{paperA}
\bibliographystyle{amsalpha}
\end{document}